\documentclass[11pt]{article}

\usepackage{amsmath}
\usepackage{amsopn}
\usepackage{latexsym}
\usepackage{amsfonts}
\usepackage{amssymb}
\usepackage{graphicx}
\usepackage{algorithm}
\usepackage{algorithmic}
\usepackage{float}
\usepackage{xcolor}

\newcommand{\eps}{\varepsilon}

\DeclareMathOperator*{\argmin}{\mbox{arg}\min}
\DeclareMathOperator*{\argmax}{\mbox{arg}\max}

\newtheorem{proposition}{Proposition}

\newtheorem{lemma}{Lemma}

\begin{document}

\title{On Median Filters\\ for Motion by Mean Curvature}
\author{Selim Esedo\=glu\thanks{Corresponding author. 530 Church St, Ann Arbor, MI 48109. {\bf esedoglu@umich.edu}}\\ University of Michigan \and Jiajia Guo\\ University of Michigan \and David Li\\ University of Michigan}
\maketitle

\begin{abstract}
  The median filter scheme is an elegant, monotone discretization of the level set formulation of motion by mean curvature.
It turns out to evolve every level set of the initial condition precisely by another class of methods known as threshold dynamics.
  Median filters are, in other words, the natural level set versions of threshold dynamics algorithms.
  Exploiting this connection, we revisit median filters in light of recent progress on the threshold dynamics method.
  In particular, we give a variational formulation of, and exhibit a Lyapunov function for, median filters, resulting in energy based unconditional stability properties.
The connection also yields analogues of median filters in the multiphase setting of mean curvature flow of networks.
These new multiphase level set methods do not require frequent redistancing, and can accommodate a wide range of surface tensions.
\end{abstract}

\section{Introduction}
\label{sec:Introduction}
\let\thefootnote\relax\footnotetext{{\it 2020 Mathematics Subject Classification.} Primary: 65M12; Secondary 35K93.}
Motion by mean curvature of an interface, or networks of interfaces (curves in $\mathbb{R}^2$, or surfaces in $\mathbb{R}^3$) arises in a great variety of applications.
Formally, it can be seen as gradient flow (steepest descent) dynamics associated with surface tension, i.e. perimeter (length in $\mathbb{R}^2$, or surface area in $\mathbb{R}^3$).
At any point on an interface, away from free boundaries known as triple junctions in the network case, the interface moves in the normal direction with speed proportional to its mean curvature at that point.
If the interface is parametrized, the evolution is described by a system of parabolic partial differential equations satisfied by the time dependent parametrization.
Prominent among the applications this dynamics plays a central role in are materials science (evolution of microstructure in polycrystalline materials, see e.g. \cite{mullins, kinderlehrer1, kinderlehrer4, kinderlehrer5}) and computer vision (variational models for image segmentation, see e.g. \cite{mumford_shah, chan_vese, vese_chan}).
In both, as well as in many others, the more challenging network (also known as the multiphase) case of the problem is relevant.
Here, there is the additional complication of enforcing natural boundary conditions at triple junctions along which three distinct interfaces intersect.

Both in the scalar two-phase and especially in the vectorial multi-phase case, the evolution can entail singularities and topological changes.
All these features make efficient and accurate simulation of the dynamics challenging.
Numerical methods that rely on explicit, parametrized representation of the interfaces face the daunting task of detecting, classifying, and ``manually'' carrying out topological changes that are all but inevitable in the long run.
Implicit methods, such as phase field, level set \cite{osher_sethian}, and threshold dynamics \cite{merriman_bence_osher}, on the other hand, have to contend with nonlinear, degenerate, or singular PDE descriptions, and have to ensure correct conditions at free boundaries are indirectly enforced.

In this paper, we focus on two related approaches that represent the interfaces implicitly: The level set method \cite{osher_sethian}, and threshold dynamics \cite{merriman_bence_osher, mbo92}.
The former has a complete theory in the scalar (two-phase) case for which elegant, convergent (with proof) numerical methods have also been developed that are, nevertheless, typically low order accurate in time and come with oppressive time step restrictions.
A major advantage is spatial accuracy on uniform grids: The interface can be located at sub-grid precision via interpolation.
Its main weakness is in extension to networks (the multi-phase case) at the generality demanded by applications, where verifying that correct junction conditions hold at free boundaries remains far from obvious.
The latter is unconditionally stable, with very low per time step cost, and has a fairly complete theory even in the network (multi-phase) case, including the verification of the correct junction conditions along free boundaries.
Its main drawback is low accuracy when implemented naively on uniform spatial grids: Because it represents interfaces by characteristic functions, the interface cannot be located at subgrid precision via interpolation.
As a result, even a smooth interface the curvature of which is small enough compared to the choice of time step size can get ``pinned'' (stuck).

Exploiting a precise connection between a particularly elegant discretization of the level set method known as the median filter scheme \cite{oberman2004,morel_book} and threshold dynamics, the present study offers the following contributions:
\begin{enumerate}
\item A new, variational formulation of median filter schemes for the level set method, in any dimension.
In particular, a Lyapunov functional is given in Section \ref{sec:Variational} that implies unconditional energy stability.
A minimizing movements interpretation is also offered.
In earlier work, the comparison principle had been the main tool for investigating stability.
A potential application is indicated.
\item A new, fast and accurate algorithm for approximating median filters in two dimensions (three dimensions viable, left to future work), in Section \ref{sec:Fast}.
In the scalar (two-phase) case, convergence of this approximation to the unique viscosity solution is verified.
\item A new, {\em monotone} discretization of the level set equation for mean curvature motion that is {\em second order accurate} in space and time in dimension two.
Convergence to the unique viscosity solution is again ensured; Section \ref{sec:HighOrder}.
\item Some barriers to finding second order accurate in time versions of median filter schemes in dimensions three and higher; Section \ref{sec:barrier}.
\item New, multiphase analogues of median filter schemes in any dimension, based entirely on exploiting the connection to threshold dynamics and our rather complete understanding of the latter in the multiphase context of networks \cite{esedoglu_otto, laux_otto, laux_otto2}; Section \ref{sec:Multiphase}.
This results in a new level set method for multiphase mean curvature motion that allows locating the interface via interpolation and enforces the correct junction condition at the free boundaries, at the generality demanded by applications (e.g. the unequal, non-additive surface tension case \cite{esedoglu_otto}).
\end{enumerate}
In addition to presenting new analysis and extensions of median filters (and hence the level set method) as summarized above, the present work can also be seen as a contribution to threshold dynamics, in addressing its difficulties on uniform grids by finding natural level set versions of it.
\section{Background}
\label{sec:Background}
\noindent The level set formulation of two-phase motion by mean curvature is described by the partial differential equation
\begin{equation}
\label{eq:levelset}
\phi_t = |\nabla\phi| \nabla \cdot \left( \frac{\nabla\phi}{|\nabla\phi|} \right)
\end{equation}
A complete well-posedness theory is developed in \cite{evans_spruck, chen_giga_goto} in the framework of viscosity solutions.
Discretization of the equation, which is degenerate parabolic, and singular wherever $\nabla \phi = 0$, has been an interesting problem of numerical analysis.
Preserving qualitative features of the viscosity solution, such as the comparison principle that forms its backbone, is an important challenge in the design of numerical schemes.

Among the many interesting contributions to the numerical treatment of (\ref{eq:levelset}), the {\em local median filter} based algorithm proposed in \cite{oberman2004} is one of the most elegant (see also \cite{morel_book} and references therein for earlier related work).
If we denote by $\mathbf{M}_r\phi^n(x)$ the median of the level set values at time step $n$ around the periphery $\partial B_r(x)$ of the ball $B_r(x)$ of radius $r$ centered at $x$, the scheme is simply
\begin{equation}
\label{eq:ob}
\phi^{n+1}(x) = \mathbf{M}_r \phi^n(x)
\end{equation}
It was derived from the level set formulation (\ref{eq:levelset}) based on the observation that the right hand side can be written, at least when $\phi$ is sufficiently differentiable and $\nabla\phi(x)\not= 0$, as the Laplacian of $\phi$ in the tangent plane of its level set passing through $x$:
\begin{equation}
|\nabla\phi| \nabla \cdot \left( \frac{\nabla\phi}{|\nabla\phi|} \right) = \Delta \phi(x) - \left\langle D^2 \phi(x) \frac{\nabla\phi}{|\nabla\phi|} \, , \, \frac{\nabla\phi}{|\nabla\phi|} \right\rangle
\end{equation}
Based on this, scheme (\ref{eq:ob}) can immediately be seen to be at least plausible in dimension $d=2$ by observing that
\begin{equation}
\phi \left( x \pm \frac{\nabla^\perp \phi(x)}{|\nabla^\perp \phi(x)|}r \right) = \mathbf{M}_r\phi(x) + O(r^3) \mbox{ as } r\to 0
\end{equation}
so that
\begin{equation}
\left\langle D^2 \phi(x) \frac{\nabla^\perp\phi(x)}{|\nabla^\perp\phi(x)|} \, , \, \frac{\nabla^\perp\phi(x)}{|\nabla^\perp\phi(x)|} \right\rangle \approx \frac{\mathbf{M}_r\phi(x) - 2\phi(x) + \mathbf{M}_r\phi(x)}{r^2}.
\end{equation}
Equation (\ref{eq:ob}) then implies
\begin{equation}
\phi^{n+1}(x) = \phi^n(x) + k |\nabla\phi^n(x)| \nabla \cdot \left( \frac{\nabla\phi^n(x)}{|\nabla\phi^n(x)|} \right) + O(k^2)
\end{equation}
provided we choose the time step size $k$ as $k=\frac{1}{2} r^2$.

A desirable and most helpful property of scheme (\ref{eq:ob}) is its {\em monotonicity}: 
\begin{equation}
\label{eq:monotonicity}
\phi_1^0(x) \geq \phi_2^0(x) \mbox{ for all } x \Longrightarrow \phi_1^n \geq \phi_2^n(x) \mbox{ for all } x \mbox{ and } n,
\end{equation}
which of course also holds for the viscosity solutions of the PDE (\ref{eq:levelset}).
Another useful property, especially for numerical implementation, is that any global Lipschitz bound is preserved:
\begin{equation}
\sup_{x \not= y} \frac{|\phi^n(x) - \phi^n(y)|}{|x-y|} \leq \sup_{x \not= y} \frac{|\phi^0(x) - \phi^0(y)|}{|x-y|}.
\end{equation}
This implies, in particular, that no {\em steepening} of the level-set function will take place, which is to be expected since {\em every} level set approximates motion by mean curvature, which is well known to enjoy a related property.

%\begin{center}
%\fbox{\begin{minipage}{4.5 in}
%{\bf Algorithm:} Fix a time step size $\delta t>0$. Let $r = \frac{1}{2} \delta t^2$.
%Generate discrete in time approximations $\phi^k(x) $ to (\ref{eq:levelset}) at times $t=k\delta t$ with the following update: 
%\begin{enumerate}
%\item Sample the function $\phi^k(x)$ uniformly along $\partial B_r(x)$: 
%  \begin{equation}
%    \label{eq:ob1}
%    \mathcal{S} = \{ \phi^k(y_j) \, : \, y_j\in \partial B_r(x) \mbox{ for } j=1,2,\ldots,m. \}
%  \end{equation}
%\item Set
%  \begin{equation}
%    \label{eq:ob2}
%    \phi^{k+1}(x) = \median(\mathcal{S})
%    \end{equation}
%\end{enumerate}
%\end{minipage}}
%\end{center}

In practice, to implement scheme (\ref{eq:ob}), we can choose $m$ points $y_1,y_2,\ldots,y_m$ that sample approximately uniformly the sphere $\partial B_r(0)$, with $r=\sqrt{2k}$.
The value of the level set function $\phi$, which is typically presented on a uniform grid, can be evaluated at any $x+y_j$ via bilinear interpolation \cite{ryo}, which preserves monotonicity.
The algorithm is then as follows:

\begin{algorithm}[H]
\caption{Median Filter for Motion by Mean Curvature \cite{oberman2004}}
\label{alg:ob}
\begin{algorithmic}[1]
\STATE {\bf sort} the level set values $\{ \phi^n(x+y_1), \phi^n(x+y_2), \ldots, \phi^n(x+y_m)\}$ so that the permutation $p \, : \, \{1,2,\ldots,m\} \to \{1,2,\ldots,m\}$ satisfies
$$ \phi^n \big( x+y_{p(1)} \big) \leq \phi^n \big( x+y_{p(2)} \big) \leq \cdots \leq \phi^n \big( x+y_{p(m)} \big).$$
\STATE $ \phi^{n+1}(x) = \frac{1}{2} \Big( \phi^n(x+y_{p(\lfloor \frac{m}{2} \rfloor)}) + \phi^n(x+y_{p(\lceil \frac{m}{2} \rceil)}) \Big)$.
\end{algorithmic}
\end{algorithm}

The main task in Algorithm \ref{alg:ob} is {\em sorting} the values $\{ \phi^n(x+y_j)\}_{j=1}^m$.
In what follows, it will be useful to regard this reliance on the sort operation (applied to level set values) as the distinctive feature of median based scemes; doing so will facilitate some of the extensions that will be introduced.

Algorithm \ref{alg:ob} is discrete in time but continuous in space, save for the discrete set of points $y_j$ sampled on $\partial B_r(0)$.
It is convenient to consider a discrete in time, but fully continuous in space version that relies on the {\em continuum} median.
In doing so, it's also worth discussing a slight generalization, namely the {\em weighted local median}.
To that end, let $K$ be a positive, radially symmetric kernel with rapid decay and unit mass.
We will informally allow $K$ to be concentrated (a $\delta$ function) on finitely many circles ($d=2$) or spherical shells ($d=3$) to cover algorithms such as Algorithm \ref{alg:ob}.
For $r>0$, define
\begin{equation}
K_r(x) = \frac{1}{r^d} K \left( \frac{x}{r} \right).
\end{equation}
%\begin{equation*}
%\phi^{k+1}(x) \in \argmin_\xi \int_{\partial B_1(0)}  \big| \xi - \phi(x+ry) \big| \, dS(y)
%\end{equation*}
%Therefore, we will consider the more general  {\em weighted local median}
%\begin{equation*}
%\mbox{\bf M}_K \phi (x) = \argmin_\xi \int \left| \xi - \phi(x+y) \right| K(y) \, dy
%\end{equation*}
%where we informally allow the radially symmetric, positive kernel $K$ to be concentrated (a delta function) on unions of circles (for $d=2$) or spherical shells (for $d=3$).
For a bounded, continuous function $\phi$ and $\lambda\in\mathbb{R}$, let
\begin{equation}
\mathbf{T}_\lambda \phi = \big\{ x \, : \, \phi(x) \geq \lambda \big\}.
\end{equation}
The function
\begin{equation}
\label{eq:psi}
\mathbf{\psi}_K \phi (x,\lambda) = \int_{T_\lambda \phi} K(x-y) \, dy
\end{equation}
is decreasing (since $K\geq 0$) and left continuous in $\lambda$, satisfying
$$ \lim_{\lambda\to\infty} \mathbf{\psi}_K \phi (x,\lambda) = 0 \mbox{ and } \lim_{\lambda\to -\infty} \mathbf{\psi}_K \phi (x,\lambda) = 1. $$
Define the weighted median of $\phi$ at $x$ with respect to the weight $K$ as
\begin{equation}
\label{eq:median}
\mathbf{M}_K\phi(x) = \sup \left\{ \lambda \, : \, \psi_K\phi(x,\lambda) \geq \frac{1}{2} \right\}.
\end{equation}
If $\phi$ is continuous, then so is $\mathbf{M}_K\phi$.
Moreover, $\mathbf{M}_K$ preserves any global Lipschitz bound, $L$,  that may hold for $\phi$:
\begin{equation}
\label{eq:lipbound}
\psi_K\phi(x,\lambda+L|x-x'|) \leq \psi_K\phi(x',\lambda) \leq \psi_K\phi(x,\lambda-L|x-x'|)
\end{equation}
for any $x,x',\lambda$, which implies
\begin{equation}
\label{eq:lipbound2}
\mathbf{M}_K \phi(x)-L|x-x'| \leq \mathbf{M}_K\phi(x') \leq \mathbf{M}_K\phi(x) + L|x-x'|.
\end{equation}
Note also the following well-known variational characterization of the weighted local median:
\begin{equation}
\label{eq:medianvariational}
\mathbf{M}_K \phi(x) \in \argmin_\xi \int \big| \xi - \phi(y) \big| K(x-y) \, dy
\end{equation}
Indeed, we have
\begin{equation}
\big| \xi - \phi(y) \big| = \int_{-\infty}^\xi \mathbf{1}_{(\mathbf{T}_\lambda\phi)^c} (y) \, d\lambda + \int_\xi^\infty \mathbf{1}_{\mathbf{T}_\lambda \phi}(y) \, d\lambda
\end{equation}
so that
\begin{equation}
\label{eq:medianvariational2}
\begin{split}
C(\xi) :&= \int \big| \xi - \phi(x) \big| K(x-y) \, dx\\
&= \int_{-\infty}^\xi \int_{\big( \mathbf{T}_\lambda \phi \big)^c} K(x-y) \, dy \, d\lambda + \int_\xi^\infty \int_{\mathbf{T}_\lambda \phi} K(x-y) \, dy \, d\lambda\\
&= \int_{-\infty}^\xi 1 - \psi_K\phi(x,\lambda) \, d\lambda + \int_\xi^\infty \psi_K\psi(x,\lambda) \, d\lambda
\end{split}
\end{equation}
Since $\psi_K\phi(x,\lambda)$ is a decreasing function of $\lambda$, with $\lim_{\lambda\to\infty} \psi_K\phi(x,\lambda) = 0$ and $\lim_{\lambda\to-\infty} \psi_K\phi(x,\lambda) = 1$, equation (\ref{eq:medianvariational2}) shows the continuous function $C(\xi)$ is decreasing on $(-\infty,\mathbf{M}_K\phi(x))$ and increasing on $(\mathbf{M}_K\phi(x),\infty)$, which implies statement (\ref{eq:medianvariational}).

The discrete in time, continuous in space version of weighted median filter scheme for motion by mean curvature is the following generalization of (\ref{eq:ob}):
\begin{equation}
\label{eq:obw}
\phi^{n+1}(x) = \mathbf{M}_{K_{\sqrt{2k}}} \phi^n(x).
\end{equation}
In the fully discrete setting, the weighted local median can be computed with a sort operation just like the standard median.
For instance, when $K$ is a continuous, radially symmetric, compactly supported positive function, let $y_j$ for $j=1,2,\ldots,m$ sample approximately uniformly the support of $K$.
Then, we may approximate $\mathbf{M}_K\phi(x)$, and hence scheme \ref{eq:obw}), as described in Algorithm \ref{alg:wmf} below:

\begin{algorithm}[H]
\caption{Weighted Median Filter \cite{morel_book}}
\label{alg:wmf}
\begin{algorithmic}[1]
\STATE {\bf sort} the level set values $\{ \phi^n(x+y_1), \phi^n(x+y_2), \ldots, \phi^n(x+y_m)\}$ so that the permutation $p \, : \, \{1,2,\ldots,m\} \to \{1,2,\ldots,m\}$ satisfies
$$ \phi^n \big( x+y_{p(1)} \big) \leq \phi^n \big( x+y_{p(2)} \big) \leq \cdots \leq \phi^n \big( x+y_{p(m)} \big).$$
\STATE {\bf set} $C=0$ and $\ell = 1$.
\WHILE{$ C < \frac{1}{2} \|K\|$}
\STATE $C \leftarrow C + K \big(y_{p(\ell)} \big)$
\STATE $\ell \leftarrow \ell+1$
\ENDWHILE
\STATE $ \phi^{n+1}(x) = \frac{1}{2} \Big( \phi^n(x+y_{p(\ell-1)}) + \phi^n(x+y_{p(\ell)}) \Big)$
\end{algorithmic}
\end{algorithm}
\noindent Steps 2 through 6 of Algorithm \ref{alg:wmf}, which is its only difference from Algorithm \ref{alg:ob}, is a simple cumulative sum of kernel values, in the order obtained from the sort operation of Step 1.
Even the multiphase median filter, namely Algorithm \ref{alg:mmf}, discussed in Section \ref{sec:Multiphase}, differs from Algorithm \ref{alg:ob} merely in the manner of this cumulative sum.

A different, also very elegant, algorithm for simulating mean curvature motion, including in the multiphase setting, is threshold dynamics \cite{mbo92, merriman_bence_osher}.
In this method, the interface, which is the boundary of a set, is represented by the characteristic function of that set.
Given the initial set $\Sigma^0$, the algorithm generates a discrete in time approximation to the flow as follows:

\begin{algorithm}[H]
\caption{Two-Phase Treshold Dynamics \cite{merriman_bence_osher,mbo92}}
\label{alg:td}
\begin{algorithmic}[1]
\STATE Convolution step:
\begin{equation}
\label{eq:td1}
u = K_{k} * \mathbf{1}_{\Sigma^n}
\end{equation}
\STATE Thresholding step:
\begin{equation}
\label{eq:td2}
\Sigma^{n+1} = \Big\{ x \, : \, u(x) \geq \frac{1}{2} \Big\}.
\end{equation}
\end{algorithmic}
\end{algorithm}

Let us introduce notation for applying one step of this scheme to a set $\Sigma$:
\begin{equation}
  \label{eq:T}
  \mathbf{S}_{K} \Sigma = \Big\{ x \, : \, K * \mathbf{1}_{\Sigma}(x) \geq \frac{1}{2} \Big\}
\end{equation}
so that Algorithm (\ref{eq:td1}) \& (\ref{eq:td2}) can be expressed as
\begin{equation}
  \Sigma^{n+1} = \mathbf{S}_{K_{k}} \Sigma^n = \mathbf{T}_{\frac{1}{2}} K_k * \mathbf{1}_{\Sigma^n}
\end{equation}
In summary, we have two distinct classes of algorithms, namely median filters and threshold dynamics, for motion by mean curvature, represented by the maps $\mathbf{M}_K$ and $\mathbf{S}_K$, respectively.

Next, we recall the variational formulation for threshold dynamics algorithm (\ref{eq:td1}) \& (\ref{eq:td2}) from \cite{esedoglu_otto}.
We have:
\begin{equation}
\label{eq:tdvariational}
\begin{split}
\Sigma^{n+1} &= \mathbf{T}_{\frac{1}{2}} K * \mathbf{1}_{\Sigma^n} = \argmin_\Sigma \int_\Sigma 1 - 2 K * \mathbf{1}_{\Sigma^n} \, dx\\
&= \argmin_\Sigma \int_{\Sigma^c} K * \mathbf{1}_{\Sigma} \, dx  + \int \Big( \mathbf{1}_\Sigma - \mathbf{1}_{\Sigma^n} \Big) K * \Big( \mathbf{1}_\Sigma - \mathbf{1}_{\Sigma^n} \Big) \, dx
\end{split}
\end{equation}
When the Fourier transform $\widehat{K}$ of the kernel $K$ is positive, the term
\begin{equation}
\label{eq:tdml}
\int \Big( \mathbf{1}_\Sigma - \mathbf{1}_{\Sigma^n} \Big) K * \Big( \mathbf{1}_\Sigma - \mathbf{1}_{\Sigma^n} \Big) \, dx
\end{equation}
is positive and achieves its minimum value of $0$ at $\Sigma = \Sigma^n$.
It thus acts as a {\em movement limiter}, i.e. threshold dynamics is thus seen to be a {\em minimizing movement}  \cite{degiorgiminmov,almgren_taylor_wang,luckhaus95} for the energy
\begin{equation}
\label{eq:E}
E_K(\Sigma) = \int_{\Sigma^c} K * \mathbf{1}_{\Sigma} \, dx
\end{equation}
and hence dissipates it at every time step:
\begin{equation}
\label{eq:EO}
E_K(\mathbf{S}_K \Sigma) \leq E_K(\Sigma).
\end{equation}

It's worth recalling from the Introduction a difficulty with threshold dynamics algorithms such as Algorithm \ref{alg:td} that is {\em not} shared by median filters: Because they represent interfaces via discontinuous (characteristic) functions, threshold dynamics schemes can be awkward on uniform grids.
When the time step size is small compared to the curvature of the interface, there can be artificial pinning (the interface gets stuck) with naive implementations.
Numerous workarounds are available, including an adaptive version \cite{ruuththesis,ruuth0} that relies on non-uniform FFT for step (\ref{eq:td2}); these are effective, but brute force solutions.
One of the goals of the present study is a step towards a more fundamental solution, perhaps via the connection with median filters, in applications such as multiphase motion by mean curvature of networks (Section \ref{sec:Multiphase}).

\section{Connection to Threshold Dynamics}
\label{sec:Connection}
In this section, we record the important observation that, in fact, threshold dynamics (Algorithm \ref{alg:td}) and weighted median filtering (Algorithm \ref{alg:wmf}) are in some precise sense equivalent.
Indeed, the median based Algorithm \ref{alg:ob}, and especially its sort-based immediate generalizations discussed in Section \ref{sec:Background}, {\em precisely carry out threshold dynamics} on every super level set $\{\phi \geq \lambda \}$ of the level set function $\phi$.
We state it clearly in the notation of this paper:

\begin{lemma}
\label{connection}
Let $K$ be positive and rapidly decaying.
Let $\phi$ be continuous and bounded.
For every $\lambda\in\mathbb{R}$, we have
\begin{equation}
\label{eq:sttm}
\mathbf{S}_{K} \mathbf{T}_\lambda \phi = \mathbf{T}_\lambda \mathbf{M}_K \phi
\end{equation}
and therefore
\begin{equation}
\label{eq:stacked}
\mathbf{M}_K\phi(x) = \sup \big\{ \lambda \, : \, x \in \mathbf{S}_K \mathbf{T}_\lambda \phi \big\}.
\end{equation}
\end{lemma}
\medskip

\noindent {\bf Proof:}
This is simple, but worth spelling out carefully.
Let $x\in\mathbf{S}_K\mathbf{T}_\lambda\phi$.
This means
\begin{equation}
  K * \mathbf{1}_{ \{ \phi \geq \lambda \} } (x)  = \int_{T_\lambda \phi} K(x-y) \, dy = \psi_K \phi(x,\lambda) \geq \frac{1}{2}.
\end{equation}
Since $\psi_K \phi (x,\xi) < \frac{1}{2}$ for all $\xi>\mathbf{M}_K\phi(x)$, we have
\begin{equation}
  \mathbf{M}_K\phi(x) \geq \lambda.
\end{equation}
Hence $x\in\mathbf{T}_\lambda \mathbf{M}_K \phi$.
Conversely, if $x\in \mathbf{T}_\lambda \mathbf{M}_K\phi$, i.e. if $\mathbf{M}_K\phi(x) \geq \lambda$, then $\psi_K \phi(x,\lambda) \geq \frac{1}{2}$ thanks to the left continuity of $\psi_K \phi(x,\cdot)$.
But then $K*\mathbf{1}_{T_\lambda\phi}(x) \geq \frac{1}{2}$, which means $x\in\mathbf{S}_K\mathbf{T}_\lambda\phi$, establishing (\ref{eq:sttm}), from which (\ref{eq:stacked}) follows.
$\Box$
\medskip

Thus, in the terminology of \cite{morel_book}, median schemes are the {\em stacked filter} versions of threshold dynamics schemes.
Level by level application of the threshold dynamics scheme to a continuous level set function appears numerous times in earlier studies focused on the analysis of these algorithms, including the pioneering works \cite{evans_mbo, ishii_pires_souganidis} that furnished some of the first convergence proofs.
In these works, the connection to a median operation on the level set function is not explicitly stated.
However, it turns out the connection was already observed in subsequent earlier work, e.g. in \cite{morel_book}.

\section{Variational Formulation of Median Schemes}
\label{sec:Variational}
The simple observation of Lemma \ref{connection} allows us to extend to median \& sorting based level set schemes the variational formulation of threshold dynamics introduced in \cite{esedoglu_otto}.
In particular, we will exhibit Lyapunov functions for, and thereby establish unconditional energy stability of, median and sort based level set schemes.
To that end, first note that
\begin{equation}
\int_{\mathbb{R}} \mathbf{1}_{\mathbf{T}_\lambda \phi} (y) \Big( 1 - \mathbf{1}_{\mathbf{T}_\lambda \phi} (x) \Big) \, d\lambda = \Big( \phi(y) - \phi(x) \Big)_+.
\end{equation}
Therefore,
\begin{equation}
\label{eq:note1}
\int_{\mathbb{R}} E_K \big( \mathbf{T}_\lambda \phi \big) \, d\lambda = \frac{1}{2} \iint K(x-y) |\phi(x) - \phi(y)| \, dx \, dy.
\end{equation}
This establishes the following Lyapunov function for the semi-discrete version (\ref{eq:obw}) of the median filter scheme given in Algorithm \ref{alg:wmf}:
\begin{proposition}
\label{prop:dissipation}
Let $\widehat{K} \geq 0$.
Then, the discrete in time, continuous in space version (\ref{eq:obw}) of the median filter scheme given in Algorithm \ref{alg:wmf} dissipates the non-local energy
\begin{equation}
\label{eq:TVK}
\mathcal{E}_K(\phi) = \iint K(x-y) \big|\phi(x) - \phi(y) \big| \, dx \, dy,
\end{equation}
i.e. $\mathcal{E}_K(\phi^{n+1}) \leq \mathcal{E}_K(\phi^n)$.
\end{proposition}
\smallskip

\noindent {\bf Proof:} 
\begin{equation}
\begin{split}
\mathcal{E}_K(\mathbf{M}_K\phi) &= \iint K(x-y) \int_{\mathbb{R}} \mathbf{1}_{\mathbf{T}_\lambda \mathbf{M}_K \phi} (y) \big( 1 - \mathbf{1}_{\mathbf{T}_\lambda \mathbf{M}_K \phi} (x) \big) \, d\lambda \, dx dy\\
& \mbox{(by (\ref{eq:TVK}), (\ref{eq:note1}) and (\ref{eq:E}))}\\
&= \int_{\mathbb{R}} \iint K(x-y) \mathbf{1}_{\mathbf{S}_K \mathbf{T}_\lambda \phi} (y) \big( 1 - \mathbf{1}_{\mathbf{S}_K \mathbf{T}_\lambda \phi} (x) \big) \, dx dy \, d\lambda\\
& \mbox{ (by (\ref{eq:sttm}) of Lemma \ref{connection}) }\\
&= \int_{\mathbb{R}} E_K \big( \mathbf{S}_K \mathbf{T}_\lambda \phi \big) \, d\lambda \hspace{1 cm} \mbox{ (by (\ref{eq:E}))}\\
&\leq \int_{\mathbb{R}} E_K \big( \mathbf{T}_\lambda \phi \big) \, d\lambda \hspace{1 cm} \mbox{ (by (\ref{eq:EO}))}\\
&= \mathcal{E}_K (\phi) \hspace{1 cm} \mbox{ (by (\ref{eq:TVK}) and (\ref{eq:note1}))}. \; \Box
\end{split}
\end{equation}

\noindent {\bf Remark 1:} More generally, the median filter dissipates
\begin{equation}
\int K(x-y) \big|\Phi(\phi(x)) - \Phi(\phi(y)) \big| \, dx \, dy
\end{equation}
for any increasing, absolutely continuous $\Phi:\mathbb{R}\to\mathbb{R}$, as a small modification of the calculation above shows. $\Box$
\smallskip

It is also interesting to write down how the minimizing movements formulation of threshold dynamics extends to median filters:
%\begin{proposition}
%Let $\phi^{n}$ be generated by the weighted median filter scheme (ref).
%Then:
%\begin{equation}
%\label{eq:linearized}
%\phi^{n+1} \in \argmin_\phi \int K(x-y) \big| \phi(x) - \phi^n(y) \big| \, dx \, dy
%\end{equation}
%\end{proposition}
%\smallskip

\begin{proposition}
Under the same conditions as Proposition \ref{prop:dissipation}, 
\begin{equation}
\label{eq:medianminmov}
\phi^{n+1} \in \argmin_\phi \mathcal{E}_K(\phi) + \mathcal{M}_n(\phi)
\end{equation}
where
\begin{multline}
\mathcal{M}_n(\phi)
= \iint K(x-y) \Big( 2\big| \phi(x) - \phi^n(y) \big|\\ - \big| \phi(x) - \phi(y) \big| - \big| \phi^n(x) - \phi^n(y) \big| \Big) \, dx \, dy
\end{multline}
acts as a movement limiter, i.e.
\begin{equation}
\label{eq:medianminmov2}
\mathcal{M}_n(\phi) \geq 0 \mbox{ for all } \phi \mbox{, and } \mathcal{M}_n(\phi^n) = 0.
\end{equation}
\end{proposition}
\smallskip

\noindent {\bf Proof:} Formula (\ref{eq:medianminmov}) follows immediately from (\ref{eq:medianvariational}) and (\ref{eq:TVK}).
For any two continuous functions $f,g \, : \mathbb{R}^d \to \mathbb{R}$, we have
\begin{equation}
\int \mathbf{1}_{T_\lambda f}(x) \; \mathbf{1}_{T_\lambda g}(y) \, d\lambda = \min \{ f(x) \, , \, g(y)\}.\\
\end{equation}
Therefore,
\begin{equation}
\begin{split}
\iint \Big( \mathbf{1}_{T_\lambda \phi} - \mathbf{1}_{T_\lambda \phi^n} \Big) & K * \Big( \mathbf{1}_{T_\lambda \phi} - \mathbf{1}_{T_\lambda \phi^n} \Big) \, dx \, d\lambda\\
= \iint K(x-y) \Bigg( &\min \Big\{ \phi(x) \, , \, \phi(y) \Big\}\\
+ &\min \Big\{ \phi^n(x) \, , \, \phi^n(y) \Big\}\\
- &\min \Big\{ \phi(x) \, , \, \phi^n(y) \Big\}\\
- &\min \Big\{ \phi^n(x) \, , \, \phi(y) \Big\} \Bigg) \, dx \, dy
\end{split}
\end{equation}
Noting $\min\{a,b\} = \frac{1}{2} \{ a + b - |a-b| \}$ allows canceling some terms in the last expression, which becomes
\begin{align}
\iint & \Big( \mathbf{1}_{T_\lambda \phi} - \mathbf{1}_{T_\lambda \phi^n} \Big)  K * \Big(  \mathbf{1}_{T_\lambda \phi} - \mathbf{1}_{T_\lambda \phi^n} \Big) \, dx \, d\lambda\\
&= \frac{1}{2} \iint K(x-y) \Bigg\{ 2\big| \phi^n(x) - \phi(y) \big|\\
& \hspace{3.5 cm} -\big| \phi(x) - \phi(y) \big| + \big| \phi^n(x) - \phi^n(y) \big| \Bigg\} \, dx \, dy\\
&= \mathcal{M}_n(\phi).
\end{align}
Condition $\widehat{K}\geq 0$ then implies (\ref{eq:medianminmov2}), since
\begin{equation}
\iint f K*f \, dx = \iint \widehat{K} \big| \widehat{f} \big|^2 \, d\xi \geq 0
\end{equation}
completing the proof. $\Box$
\medskip

%\noindent {\bf Remark 2:} To reiterate, in view of (\ref{eq:medianminmov}), median filter is a minimizing movement for energy (\ref{eq:TVK}) with the movement limiter
%\begin{equation}
%\int K(x-y) \Big( 2\big| \phi(x)-\phi^n(y) \big| - \big| \phi(x) - \phi(y) \big| - \big| %\phi^n(x) - \phi^n(y) \big| \Big) \, dx \, dy
%\end{equation}
%which is (apparently) positive (and zero when $\phi = \phi^n$) when $\widehat{K}\geq 0$. $\Box$
%\medskip

\noindent {\bf Remark 2:} A slight extension of the propositions above gives a sorting based, median filter type algorithm for decreasing nonlocal variational models of the form
\begin{equation}
\label{eq:NL}
NL(\phi)  = \iint K(x-y) \Big| \phi(x) - \phi(y) \Big| \, dx \, dy + \frac{\gamma}{2} \int \big( \phi(x) - f(x) \big)^2 \, dx,
\end{equation}
provided that $\widehat{K}\geq 0$.
These types of models appear in image processing, where $f \geq 0$ is an image to denoise, and $\gamma >  0$ is the user supplied {\em fidelity} parameter, see e.g. \cite{buades} (and the local model \cite{rof} that started it all).
A short argument shows that the minimizer of (\ref{eq:NL}) needs to be positive; hence we may restrict attention to $\phi\geq 0$.
Then, recall from e.g. \cite{chambolle_mcm} that the {\em fidelity term} in (\ref{eq:NL}) can be written as
\begin{equation}
\frac{1}{2} \int \big( \phi(x) - f(x) \big)^2 \, dx = \int_0^\infty \int_{\mathbf{T}_\lambda \phi} \lambda - f(x) \, dx \, d\lambda + C,
\end{equation}
where $C$ is independent of $\phi$.
Thus, (\ref{eq:NL}) can be expressed as
\begin{equation}
\label{eq:NLlayers}
NL(\phi) = \int_0^\infty \left\{ E_K(\mathbf{T}_\lambda \phi) + \gamma \int_{\mathbf{T}_\lambda \phi} \lambda - f(x) \, dx  \right\} \, d\lambda + C.
\end{equation}
\medskip
A threshold dynamics scheme to dissipate the geometric problem in the integrand of (\ref{eq:NLlayers}), namely
\begin{equation}
\label{eq:geometric}
\min_\Sigma E_K(\Sigma) + \gamma \int_{\Sigma} \lambda - f(x) \, dx
\end{equation}
can be derived by following the strategy in \cite{esedoglu_otto}: Consider the equivalent, relaxed version of (\ref{eq:geometric}):
\begin{equation}
\label{eq:relaxed}
\min_{0 \leq u \leq 1} \int (1-u) K*u + \gamma \, u \, \Big( \lambda - f(x) \Big) \, dx
\end{equation}
The condition $\widehat{K}\geq 0$ implies (\ref{eq:relaxed}) is concave, while the box constraint on $u$  is convex, meaning that minimizing the linearization (about $u = \mathbf{1}_{\mathbf{T}_\lambda \phi^n}$) of (\ref{eq:relaxed}), namely
\begin{equation}
\label{eq:linearized}
\min_{0\leq u\leq 1} \int u \Big\{ 1-2K*\mathbf{1}_{\mathbf{T}_\lambda \phi^n} + \gamma \big( \lambda - f(x) \big) \Big\} \, dx
\end{equation}
also dissipates (\ref{eq:relaxed}) itself.
The solution of the linearized problem (\ref{eq:linearized}) is given by thresholding:
\begin{equation}
u = 
\begin{cases}
\label{eq:NLthresh}
1 & \mbox{ if } K*\mathbf{1}_{\mathbf{T}_\lambda \phi^n} \geq \frac{1}{2} \big\{ 1 + \gamma \left( \lambda - f(x) \right) \big\} \\
0 & \mbox{ otherwise.}
\end{cases}
\end{equation}
which can be applied at all levels $\lambda$ simultaneously by a sort based procedure (i.e. the ``stacked filter'' version of (\ref{eq:NLthresh})), much as in Algorithm \ref{alg:wmf}. $\Box$
\medskip

\noindent {\bf Remark 3:}
Formula (\ref{eq:medianvariational}) already exhibits weighted median filter as a plausible optimization strategy for the Lyapunov function (\ref{eq:TVK}).
Indeed, in a fully discrete setting, the update
\begin{equation}
\phi(x) \to \argmin_\xi \int \big| \xi - \phi(y) \big| K(x-y) \, dy
\end{equation}
when applied sequentially in $x$ (i.e. in a {\em Gauss-Seidel} fashion), is easily seen to decrease energy (\ref{eq:TVK}) with every update, without a condition on the kernel $K$, as this is nothing more than coordinate descent on (\ref{eq:TVK}); this is the algorithm in \cite{li_osher} for energies such as (\ref{eq:NL}).
The content of Proposition \ref{prop:dissipation} (and Remark 2) is that the same update applied concurrently at all $x$ (i.e. in a {\em Jacobi} fashion), which is precisely the median filter (\ref{eq:obw}), also decreases the energy.
This is not as immediate; in particular, it may be interesting that the Fourier transform $\widehat{K}$ of the kernel $K$ shows up in the $L^1$ context of (\ref{eq:medianvariational}). 
Counterexamples to energy decrease when this condition is violated can be found in \cite{esedoglu_jacobs} in the context of threshold dynamics.
$\Box$
\medskip

Of course, the condition $\widehat{K}\geq 0$ is not satisfied by kernels $K$ concentrated on finitely many circles or spheres, as in the original Algorithm \ref{alg:ob}, which is the case for which we describe a fast and accurate approximation to the weighted median in Section \ref{sec:Fast}.
That approximation can in fact be extended to any smooth, compactly supported, radially symmetric kernel $K$ at a modest increase in computatioal cost.
Alternatively, we can also recall slightly ``more Gauss-Seidel'' versions of threshold dynamics algorithm (\ref{eq:td1}) \& (\ref{eq:td2}) from \cite{esedoglu_jacobs} that only require $K\geq 0$ for energy dissipation.
The corresponding median filter scheme is of the form:

\begin{algorithm}[h]
\caption{Two-Step Median Filter}
\begin{algorithmic}[1]
\STATE Grow super level-sets:
        \begin{equation}
          \label{eq:gs1}
          \phi^{n+\frac{1}{2}} = \max \big\{ \phi^n \, , \, \mathbf{M}_K\phi^n \big\}.
        \end{equation}
\STATE Shrink super level-sets:
        \begin{equation}
          \label{eq:gs2}
          \phi^{n+1} = \min \big\{ \phi^{n+\frac{1}{2}} \, , \, \mathbf{M}_K\phi^{n+\frac{1}{2}} \big\}.
        \end{equation}
\end{algorithmic}
\end{algorithm}

\noindent Variational formulation of (\ref{eq:gs1}) \& (\ref{eq:gs2}) is
\begin{equation}
\begin{split}
\phi^{n+\frac{1}{2}}(x) &= \argmin_{\phi\geq\phi^n} \int K(x-y) \big| \phi(x) - \phi^n (y) \big| \, dy,\\
\phi^{n+1}(x) &= \argmin_{\phi \leq \phi^{n+\frac{1}{2}}} \int K(x-y) \big| \phi(x) - \phi^{n+\frac{1}{2}}(y) \big| \, dy,\\
\end{split}
\end{equation}
We then have the following stability property that applies to a wide variety of weights $K$, including that of Algorithm \ref{alg:ob} from \cite{oberman2004}.
\begin{proposition}
\label{prop:gs}
Let $K$ be positive, radially symmetric, with sufficient decay.
Algorithm (\ref{eq:gs1}) \& (\ref{eq:gs2}) dissipates energy (\ref{eq:TVK}).
\end{proposition}

\section{Fast and Accurate Median Computation}
\label{sec:Fast}
This section presents an alternative to the discrete approximation for the local median used in Algorithms \ref{alg:ob} and \ref{alg:wmf} which is especially useful in implementing the second order schemes in Section \ref{sec:HighOrder}.
The natural approach of sampling the support of the kernel and finding the (weighted) median of the discrete set of values thus obtained has several drawbacks.
For one, the number of points sampled would need to be scaled up even as spatial and temporal discretization is refined just to ensure consistency.
Can the number of operations per time step be kept essentially constant, independent of the time step size?
Here, we describe an approach that ensures a mild (logarithmic) growth of computational cost with respect to the accuracy desired (as it gets small), at least in two dimensions.
Convergence to viscosity solution in the two-phase setting is preserved.

In this section, we simply write $\psi_r(x,\lambda)$ to denote the function $\psi_{K_r} \phi(x,\lambda)$ defined in (\ref{eq:psi}), with the unit mass kernel $K_r$ understood to concentrate along the circle $\partial B_r(0)$.
We also write $\mathbf{M}_r$ for $\mathbf{M}_{K_r}$.
Let $\eps>0$ be a desired error tolerance in approximating $\psi_r(x,\lambda)$.
For $0 \leq \theta_1 < \theta_2 \leq 2\pi$ with $\theta_2-\theta_1 \leq \pi/2$  , define
\begin{equation}
\label{eq:bisect1}
I_{x,r}^\eps(\theta_1,\theta_2) = 
\begin{cases}
\theta_2 - \theta_1 & \mbox{if } \theta_2-\theta_1\geq \eps \; \& \; m\geq \lambda\\
I_{x,r}^\eps(\theta_1,\gamma) + I_{x,r}^\eps(\gamma,\theta_2) & \mbox{if } \theta_2-\theta_1\geq \eps \; \& \; m<\lambda<M,\\
0 & \mbox{otherwise}
\end{cases}
\end{equation}
where
\begin{equation}
\label{eq:bisect2}
\begin{split}
\gamma &= \frac{\theta_1 + \theta_2}{2}\\
m &= \min \{ \phi(x+r(\cos\theta_1,\sin\theta_1)) \, , \, \phi(x+r(\cos\theta_2,\sin\theta_2)) \}\\
M &= \max \{ \phi(x+r(\cos\theta_1,\sin\theta_1)) \, , \, \phi(x+r(\cos\theta_2,\sin\theta_2)) \}.
\end{split}
\end{equation}
In words, $I_{x,r}^\eps(\theta_1,\theta_2)$ is thus an approximation to
\begin{equation}
\int_{\Gamma_{\theta_1,\theta_2}} \mathbf{T}_\lambda \phi(y) \, ds(y)
\end{equation}
where the circular arc $\Gamma_{\theta_1,\theta_2} = \{ y \, : \, y = x + r (\cos\theta,\sin\theta) \mbox{ for } \theta_1 \leq \theta \leq \theta_2 \}$, computed using a simple bisection procedure.
Number of operations required to compute $I_{x,r}(\theta_1,\theta_2)$ is then $O(-\log\eps)$ as $\eps\to 0$.
Approximate $\psi_r(x,\lambda)$ by
\begin{equation}
\label{eq:approxpsi}
\widetilde{\psi}_{r,\eps}(x,\lambda) = \sum_{j=0}^3 I_{x,r}^\eps \left( \frac{j\pi}{2} , \frac{(j+1)\pi}{2} \right)
\end{equation}
The function $\lambda \to \widetilde{\psi}_{r,\eps}(x,\lambda)$ is decreasing.
An approximation to the weighted local median $\mathbf{M}_{K_r}\phi(x)$, defined in (\ref{eq:median}), can be found by the bisection method on this function:
\begin{equation}
\widetilde{\mathbf{M}}_{r,\eps} \phi(x) = \sup \left\{ \lambda \, : \, \widetilde{\psi}_{r,\eps}(x,\lambda) \geq \frac{1}{2} \right\}
\end{equation}
The resulting scheme is
\begin{equation}
\label{eq:approxob}
\phi^{n+1} = \widetilde{\mathbf{M}}_{r,\eps} \mathbf{B} \phi^n
\end{equation}
where $\mathbf{B}$ denotes bilinear interpolation.
The total cost becomes $O(-\log^2\eps)$, as long as the radius $r$ of the circle $\partial B_r(0)$ is merely bounded as $\eps\to 0$.
Moreover, approximation (\ref{eq:approxpsi}) and therefore the result of this bisection procedure on $\lambda$, and hence the approximate median thus obtained, is monotone in the function $\phi$.

We next show that a version of scheme (\ref{eq:approxob}) converges to the unique viscosity solution of (\ref{eq:levelset}), as is done in \cite{oberman2004} following the strategy in \cite{barles_souganidis}.
The two essential ingredients of this strategy are consistency and monotonicity, the latter of which we have already verified above.
The simplified (for brevity) version we consider is:
\begin{equation}
\phi^{n+1} = \widetilde{\mathbf{M}}_r \phi^n \mbox{ where } \widetilde{\mathbf{M}}_r \phi(x) = \lim_{\eps\to 0^+} \widetilde{\mathbf{M}}_{r,\eps} \phi(x).
\end{equation}
Note that the limit exists, as the dependence on $\eps$ is monotone.
Let us recall the notion of consistency that is required in the context of equation (\ref{eq:levelset}) from \cite{barles_georgelin}, same as in \cite{oberman2004}:
\medskip

\noindent {\bf Definition:} A numerical scheme for (\ref{eq:levelset}) of the form 
\begin{equation}
\label{eq:def1}
\phi^{n+1} = \mathcal{F}_k \phi^n \; \mbox{ i.e. } \; \frac{\phi^{n+1}-\phi^n}{k} = \frac{\mathcal{F}_k \phi^n - \phi^n}{k}
\end{equation}
is {\em consistent} if for every smooth function $\phi(x)$ it satisfies the following conditions:
\begin{enumerate}
\item At any point $x_*$ where $\nabla\phi(x_*) \not= 0$,
\begin{equation}
\label{eq:consistency1}
\lim_{k\to 0} \frac{\mathcal{F}_k \phi - \phi}{k} = |\nabla\phi| \nabla \cdot \left( \frac{\nabla \phi}{|\nabla \phi|} \right)
\end{equation}
\item At any point $x_*$ where $\nabla\phi(x_*)= 0$, let $\Lambda_1 \leq \Lambda_2$ be the eigenvalues of the Hessian $D^2\phi(x_*)$ of $\phi$.
Then,
\begin{equation}
\label{eq:consistency2}
\liminf_{k\to 0} \frac{\mathcal{F}_k \phi - \phi}{k} \geq \Lambda_1 \mbox{ and }
\limsup_{k\to 0} \frac{\mathcal{F}_k \phi - \phi}{k} \leq \Lambda_2.
\end{equation}
\end{enumerate}

\begin{figure}[h]
\begin{center}
\includegraphics[scale=0.4]{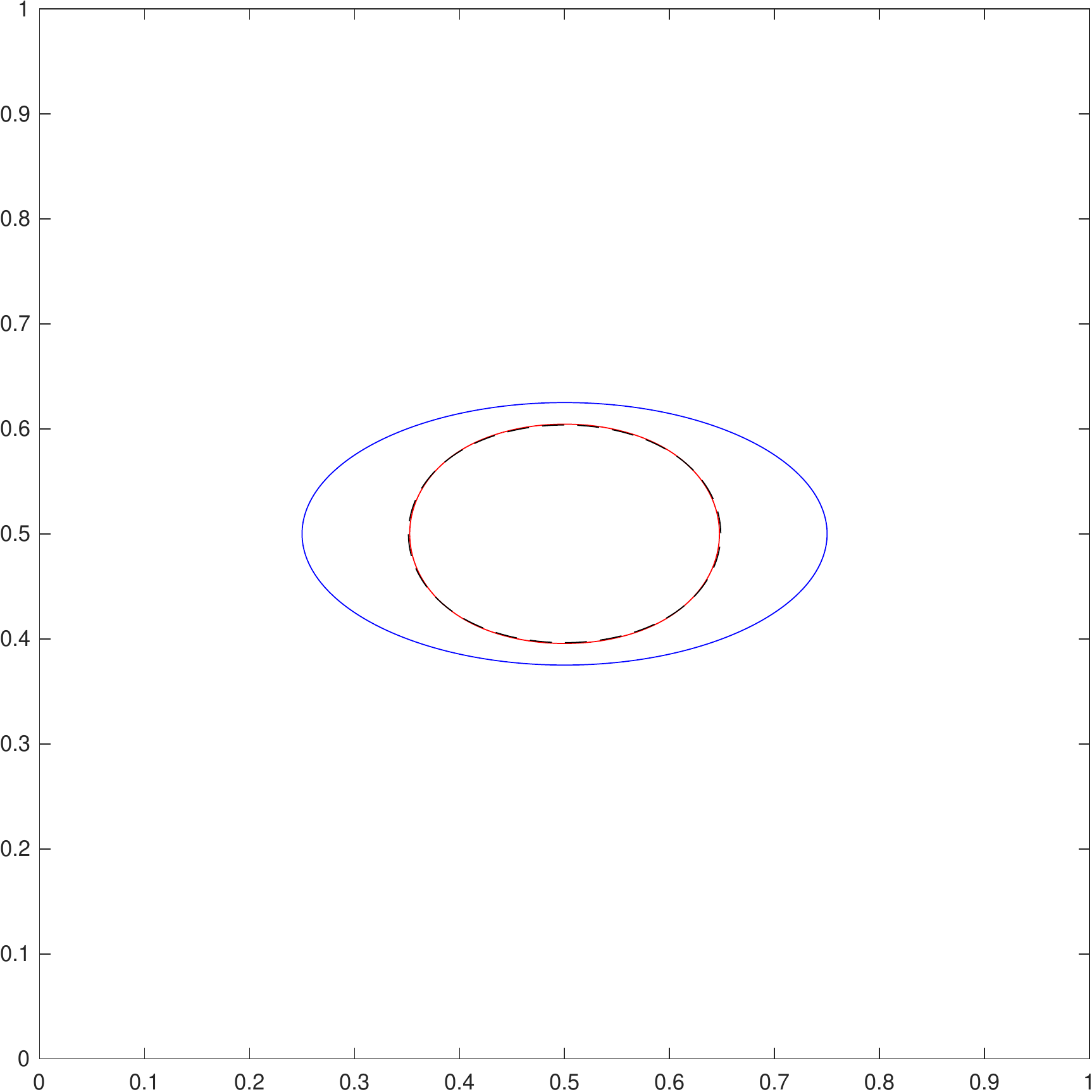}
\caption{\footnotesize Fast and accurate implementation of median filter described in Section \ref{sec:Fast}, used in Algorithm \ref{alg:ob}.
The initial curve, an ellipse, is in blue.
The final computed curve is in red.
Black is the benchmark solution, computed using front tracking with a very fine discretization of the curve.}
\label{fig:ellipse}
\end{center}
\end{figure}

\begin{table}
\begin{center}
\begin{tabular}{|c|l|c|}
\hline
$nt$ & Abs. error & Order\\
\hline
$5$ & $0.0029$ & --\\
\hline
$10$ & $0.0012$ & $1.27$\\
\hline
$20$ & $5.4\times 10^{-4}$ & $1.15$\\
\hline
$40$ & $2.51\times 10^{-4}$ & $1.1$\\
\hline
$80$ & $1.62\times 10^{-4}$ & $0.63$\\
\hline
\end{tabular}
\end{center}
\caption{\footnotesize Convergence study on a $1000\times 1000$ grid using Algorithm \ref{alg:ob} with the accurate approximation to the median described in Section \ref{sec:Fast}, exhibiting the expected first order convergence in time; see Figure \ref{fig:ellipse} for the initial and final curves.
The error is measured by an approximation to the Hausdorff distance.}
\label{tab:ellipse}
\end{table}

\begin{proposition}
\label{conv1}
Discrete in time solutions generated by the scheme
$$ \phi^{n+1} = \widetilde{\mathbf{M}}_{\sqrt{2k}} \phi^n $$
converge to the unique viscosity solution of (\ref{eq:levelset}).
\end{proposition}
\medskip

\begin{figure}
\begin{center}
\includegraphics[scale=0.5]{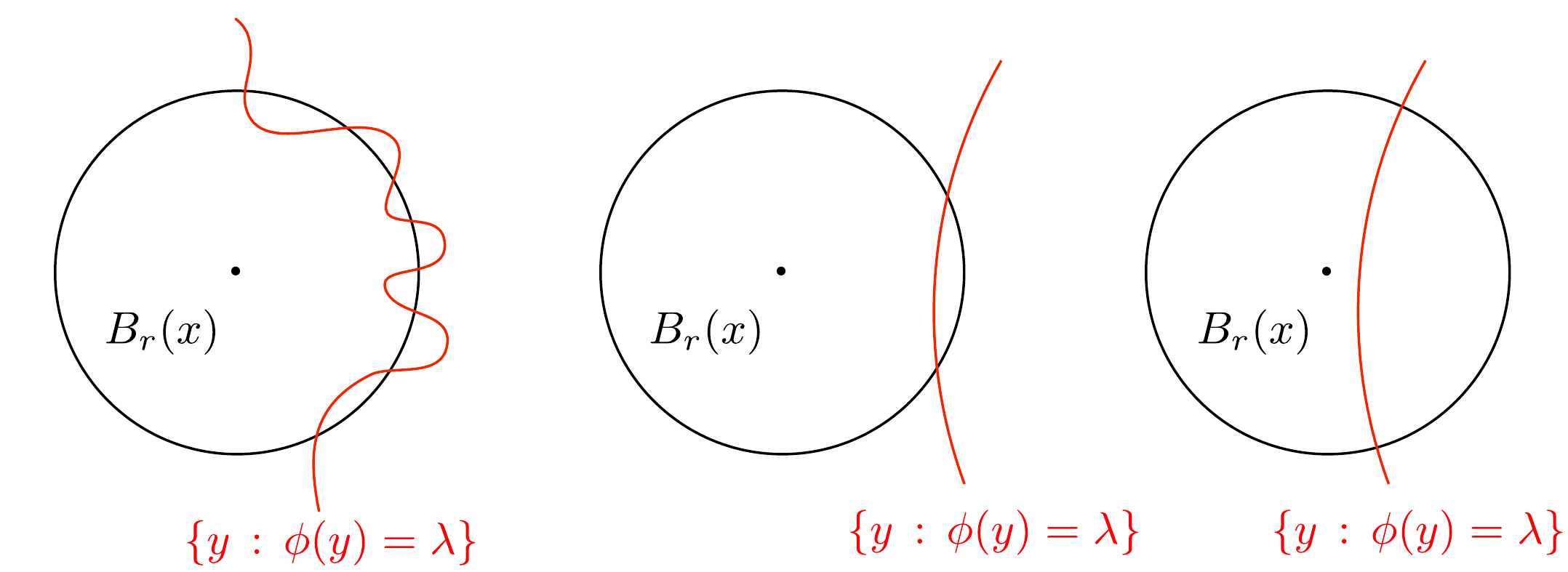}
\caption{\footnotesize The bisection procedure used in computing $\widetilde{\psi}_{r,\eps}(x,\lambda)$ can grossly miscalculate $\psi_r (x,\lambda) =  \mbox{Length}( \{ \psi = \lambda \} \cap B_r(x) )$, even at a regular point $x$ of $\phi$, regardless of $\eps>0$ for the configurations depicted on the left and in the center, but not for the configuration on the right.
Only the configuration on the right is relevant, and so convergence to viscosity solution still takes place.}
\label{fig:threeballs}
\end{center}
\end{figure}

\noindent {\bf Proof:}
We take $\mathcal{F}_k = \widetilde{\mathbf{M}}_{\sqrt{2k}}$ and verify (\ref{eq:consistency1}) \& (\ref{eq:consistency2}).
Together with monotonicity of (\ref{eq:approxob}), which has already been noted, these imply the claimed convergence via the theory in \cite{barles_souganidis}.
Let $\phi(x)$ be a smooth function.
Starting with (\ref{eq:consistency1}),
even at regular points of $\phi$, $\widetilde{\psi}_{r,\eps}(x,\lambda)$ can differ from $\psi_r(x,\lambda)$ substantially, but only if $\lambda$ is far from $\mathbf{M}_r \phi(x)$; see Figure \ref{fig:threeballs}.
Indeed, when $\nabla \phi(x_*) \not= 0$, for every $r$ small enough and $\lambda$ close enough to $\mathbf{M}_r \phi(x_*)$, the set $\{ x \, : \, \phi(x) = \lambda \}$ is a smooth curve in $B_{2r}(x_*)$ that intersects $\partial B_r(x_*)$ transversally in exactly two points that are $2r + o(r)$ apart.
This implies $\lim_{\eps\to 0} \widetilde{\psi}_{r,\eps}(x,\lambda) = \psi_r(x,\lambda)$; in words, the bisection procedure outlined in (\ref{eq:bisect1}) \& (\ref{eq:bisect2}) succeeds. 
Since $\nabla\phi(x_*) \not= 0$, that in turn implies $\widetilde{\mathbf{M}}_r \phi(x_*)= \mathbf{M}_r \phi(x_*)$.
The calculation in \cite{oberman2004}, which is summarized in Section \ref{sec:Background}, ensures consistency of $\mathbf{M}_{\sqrt{2k}} \phi(x_*)$, the exact median of $\phi$ over $\partial B_{\sqrt{2k}}(x_*)$, as an update rule for the PDE (\ref{eq:levelset}) at a regular point of a smooth function; the same therefore goes also for $\widetilde{\mathbf{M}}_{\sqrt{2k}} \phi(x_*)$, verifying (\ref{eq:consistency1}).

Turning to the case of critical points, let $\nabla\phi(x_*) = 0$.
Let $\Lambda_1 \leq \Lambda_2$ be the eigenvalues of $D^2\phi(x_*)$.
This means for any $\delta>0$, there exists $r_*(\delta) > 0$ such that if $0<r<r_*(\delta)$,
$$ \phi(x_*) + \frac{1}{2} (\Lambda_1-\delta) r^2 \leq \phi(y) \leq \phi(x_*) + \frac{1}{2} (\Lambda_2+\delta) r^2 $$
for all $y\in \partial B_r(x)$.
If $\lambda > \max_{\partial B_r(x)} \phi$, then $\widetilde{\psi}(x,\lambda) = 0$.
If $\lambda < \min_{\partial B_r(x)} \phi$, then $\widetilde{\psi}(x,\lambda) = 1$.
Thus
$$ \Lambda_1 - \delta \leq \frac{\widetilde{\mathbf{M}}_{\sqrt{2k}} \phi(x_*) - \phi(x_*)}{k} \leq \Lambda_2 - \delta, $$
which verifies (\ref{eq:consistency2}). $\Box$
\medskip

\begin{figure}[h]
\begin{center}
\includegraphics[scale=0.6]{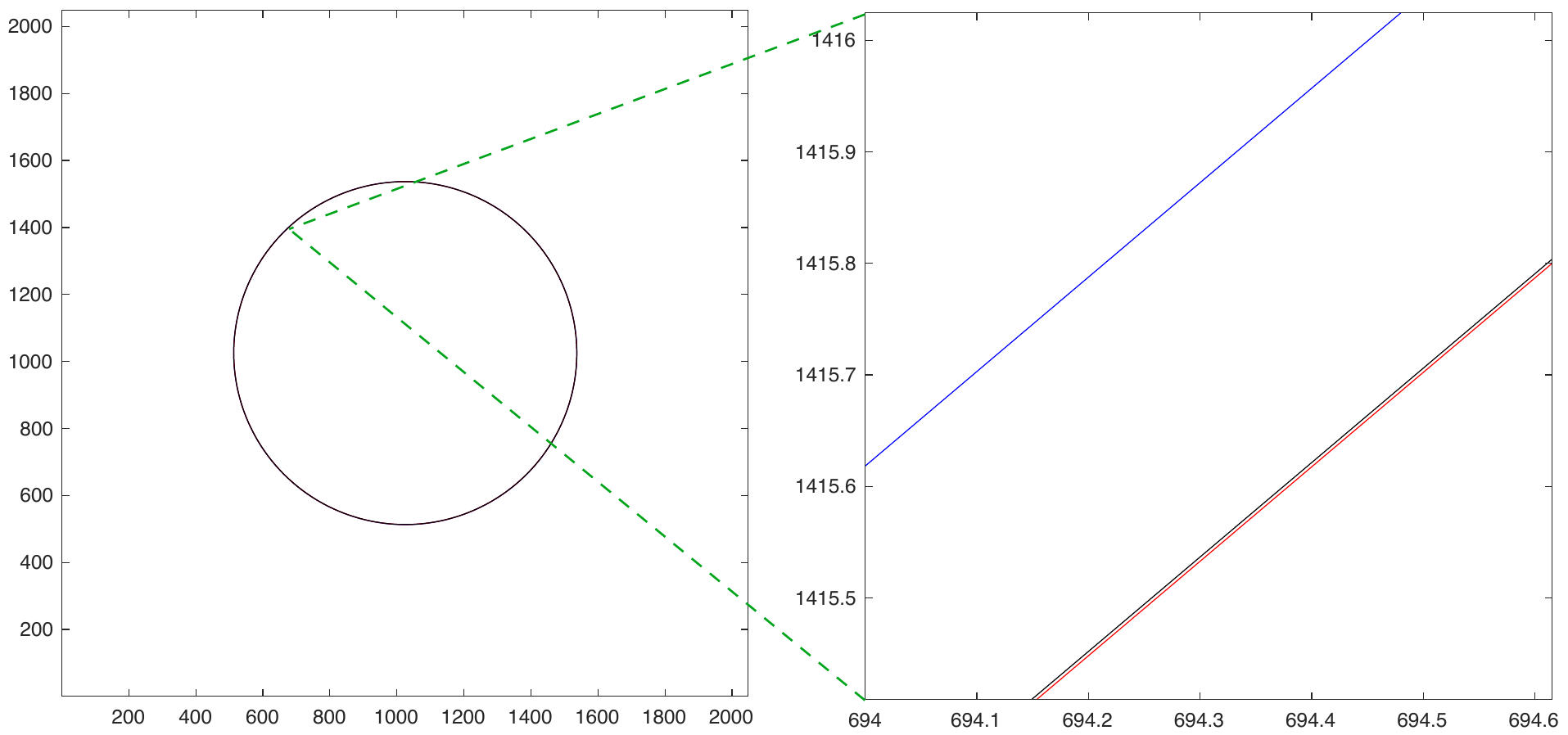}
\caption{\footnotesize Pinning of interfaces by threshold dynamics appears to be alleviated by working with their stacked filter versions, namely median filters, and using the accurate approximation of median from Section \ref{sec:Fast}.
Here, the initial curve is a very well resolved circle on a $2048 \times 2048$ grid, shown in blue.
The time step size is so small that the radius of the circular kernel in Algorithm \ref{alg:ob} is about $5$ grid points -- much smaller than naive implementations of threshold dynamics can handle without pinning.
The computed and exact curves are shown in red and black, respectively.
Perhaps because Algorithm \ref{alg:ob} is exact on circles, they are almost indistinguishable (well within $\frac{1}{100}$-th of a grid cell), even after the zoom on the right.}
\label{fig:circle}
\end{center}
\end{figure}

Figure \ref{fig:ellipse} shows a numerical experiment using Algorithm \ref{alg:ob} with the approximate median described in this section.
The initial condition is an ellipse, on a $1000\times 1000$ domain.
Table \ref{tab:ellipse} contains a convergence study as the time step size is refined, with good evidence of the expected first order convergence, until errors (measured by an approximation to the Hausdorff distance) become quite small.
As the exact solution for an initial ellipse is not explicit, the benchmark is a front tracking computation with a very fine discretization of the parametrized curve.
An even more natural numerical test would have been the exact solution of a circle; however, a short calculation shows that the discrete in time, continuous in space version of the original median filter scheme Algorithm \ref{alg:ob} is in fact {\em exact} on circles: the error is zero, regardless of the time step size.
Therefore, misleadingly high orders of convergence are observed in this test; see Figure \ref{fig:circle}.
On the other hand, this test does serve as a stark demonstration of how the pinning issue of threshold dynamics algorithms can be alleviated by working with median filters, which are their stacked versions, as discussed in Section \ref{sec:Connection}.

%%%%%%%%%%%%%%%%%%%%%%%%%%%%%%%%%%%%%%%%%%%%%%%%%%%%%%
\section{A High-Order, Monotone Median Scheme}
\label{sec:HighOrder}
A curious recent development in threshold dynamics algorithms is the second order in time, monotone version described in \cite{esedoglu_guo}.
It differs from the standard scheme only in its choice of a convolution kernel $K$:
A linear combination of three Gaussians of distinct variances is chosen that results in a positive kernel $K$, achieving second order consistency for $d=2$ (order reduces to first in $d=3$).

In this section, we explore whether it has an analogue (a monotone scheme that is second order consistent) in the realm of median filters.
Specifically, we consider a linear combination of three kernels, each concentrated on a circle of different radius, to ensure the resulting weighted local median filter can be efficiently and accurately evaluated via the method of Section \ref{sec:Fast}.
The weights in the linear combination of Gaussians in \cite{esedoglu_guo} are not all positive, even though the resulting kernel is.
Here, in forming a linear combination using kernels with disjoint supports, monotonicity clearly requires the weight associated with each kernel (circle) to be positive, so a slightly different investigation than that of \cite{esedoglu_guo} is needed.
In particular, the Appendix contains consistency calculation for convolution generated motion of an interface in $d\in\{2,3\}$ using a kernel concentrated on a circle (for $d=2$) and a spherical shell (for $d=3$).
The former calculation is used in this section to exhibit a median filter (a monotone scheme for the level set equation (\ref{eq:levelset})) that achieves second order consistency in $d=2$.
The latter is used in Section \ref{sec:barrier} that reports some bad news.

Consider a convolution kernel $K$ given as a linear combination of delta functions of various weights $c_j$ concentrated along concentric circles of radii $\alpha_j r$:
\begin{equation*}
K = \sum_j c_j \delta_{\partial B_{\alpha_j}(0)}
\end{equation*}
Using (\ref{eq:thetatot}) and solving for $\mathbf{h}$ in
\begin{equation}
\label{eq:conv2d}
\Big(\mathbf{1}_{\{(x,y) \, : \, y \geq f(x) \}} * K_r \Big) (0,\mathbf{h} r^2) = \pi \Gamma(1),
\end{equation}
where
\begin{equation}
\Gamma(p) = \sum_j c_j \alpha_j^p,
\end{equation}
using the ansatz $\mathbf{h} = \beta_0 + \beta_1 r + \beta_2 r^2 + \cdots$ we get
\begin{equation}
\begin{split}
\beta_0 =& \frac{\Gamma(2)}{2\Gamma(0)} f''(0)\\
\beta_1 =& 0\\
\beta_2 =& \left( -\frac{\Gamma^3(2) \Gamma(-2)}{48\Gamma^4(0)} + \frac{\Gamma^2(2)}{8\Gamma^2(0)} - \frac{5\Gamma(4)}{48\Gamma(0)} \right) \big(f''(0)\big)^3\\
&+ \frac{\Gamma(3)}{24\Gamma(-1)} f^{(iv)}(0).
\end{split}
\end{equation}
Choosing a linear combination of the form
$$ K = \delta_{\partial B_1(0)} + c_2 \delta_{\partial B_{\alpha_2}(0)} + c_3 \delta_{\partial B_{\alpha_3}(0)} $$
we can solve for the coefficients $c_2,\alpha_2,c_3,\alpha_3$ so that
\begin{equation}
\begin{split}
\frac{\Gamma(2)}{2\Gamma(0)} &= \frac{1}{2}\\
-\frac{\Gamma^3(2) \Gamma(-2)}{48\Gamma^4(0)} + \frac{\Gamma^2(2)}{8\Gamma^2(0)} - \frac{5\Gamma(4)}{48\Gamma(0)} &= -\frac{1}{4}\\
\frac{\Gamma(3)}{24\Gamma(-1)} &= \frac{1}{8}.
\end{split}
\end{equation}
One such solution is
$$ \alpha_2 = 2, c_2 = \frac{8}{5}, \alpha_3 = \frac{1}{2}, c_3 = \frac{32}{5} $$
resulting in the kernel
\begin{equation}
\label{eq:2ndorderkernel}
K = \delta_{\partial B_1(0)} + \frac{8}{5} \delta_{\partial B_2(0)} + \frac{32}{5} \delta _{\partial B_{\frac{1}{2}}(0)}
\end{equation}
that is positive, and gives second order consistency with the exact evolution (\ref{eq:exacttaylor}) after the redefinition $t\to \frac{1}{2}t$ of the time step size $t$.
With that, we can state the following analogue of Proposition \ref{conv1} the proof of which is also the same:

\begin{proposition}
\label{conv2}
Let the kernel $K$ be given by (\ref{eq:2ndorderkernel}).
The discrete in time scheme
\begin{equation}
\label{eq:conv2_1}
\phi^{n+1} = \widetilde{\mathbf{M}}_{K_{\sqrt{2k}}} \phi^n
\end{equation}
is monotone, and second order consistent in time.
As $k\to 0$, the solutions generated by this scheme converge to the unique viscosity solution of (\ref{eq:levelset}).
\end{proposition}
\medskip

\noindent {\bf Remark 4:} In a practical, fully discrete implementation of scheme (\ref{eq:conv2_1}), such as
\begin{equation}
\label{eq:practical}
\phi_{i,j}^{n+1} = \widetilde{\mathbf{M}}_{K_{\sqrt{2k}},\eps} \mathbf{B} \phi^n
\end{equation}
where $\mathbf{B}$ again denotes bilinear interpolation, a global truncation error of the form $O(h^2) + O(k^2)$ suggests that the natural way to scale the time step size $k$ with respect to spatial grid size $h$ is $k = O(h)$ as $h\to 0$.
Since $r=O(\sqrt{k})=O(\sqrt{h})$, the discrete median procedure of Algorithm \ref{alg:wmf} based on sampling the support of the kernel (\ref{eq:2ndorderkernel}) would entail a per time step cost of $O(N^\frac{3}{2})$ , where $N$ is the total number of grid points, i.e. $N = O(h^{-2})$.
On the other hand, taking $\eps = O(h^2)$ and using the bisection based algorithm of Section \ref{sec:Fast} on the kernel (\ref{eq:2ndorderkernel}) developed in this section, results in a fully discrete, monotone scheme that is second order consistent in both space and time, with a per time step cost of $O(N\log^2 N)$.
\smallskip

\begin{figure}[H]      
\includegraphics[scale=0.33]{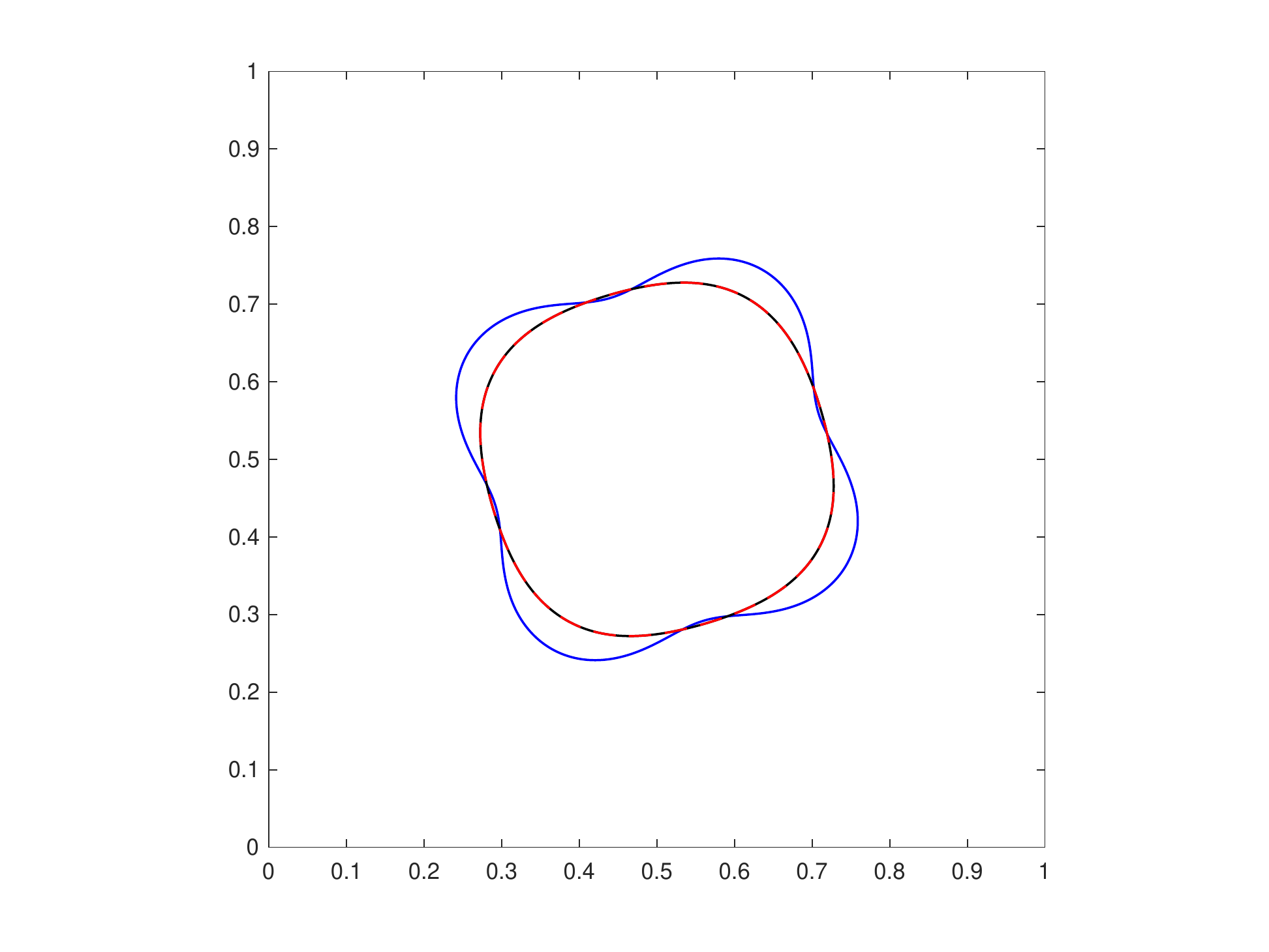} 
    \hspace{0px}
   \includegraphics[scale=0.33]{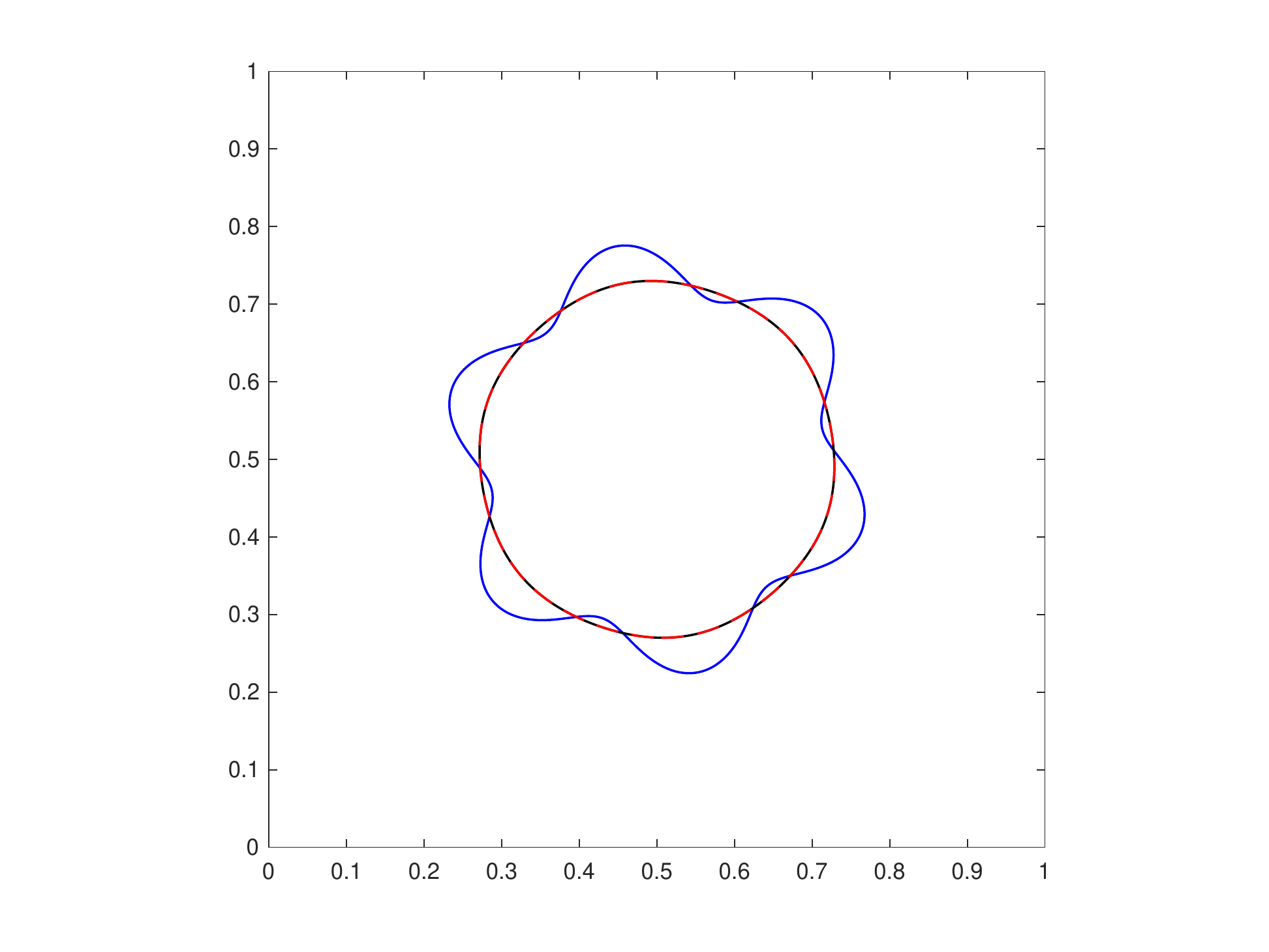}
   \caption{\footnotesize Fast and accurate implementation of weighted median filter, using the kernel (\ref{eq:2ndorderkernel}) that gives second order accuracy in time.
The initial curves, a four-petal flower and a six-petal flower, are in blue. The final computed curves are in red. Black lines are the benchmark solution, computed using front tracking with a very fine discretization of the curve.
See Table \ref{table:flowers} for the errors.}
    \label{fig:flowers}
\end{figure}

Figure \ref{fig:flowers} shows numerical experiments using kernel (\ref{eq:2ndorderkernel}) in Algorithm \ref{alg:wmf}, and an immediate extension of the fast and accurate median computation described in Section \ref{sec:Fast} to a union of three circles (the support of this new kernel).
The initial conditions are a four-petal flower and a six-petal flower.
Final time is $T = 1/200$.
The corresponding error table, Table \ref{table:flowers}, shows the measured order of convergence as the time step size and the resolution are refined simultaneously, up to a very fine $3200\times 3200$ spatial discretization, in view of the expected second order accuracy in both space and time of the practical implementation (\ref{eq:practical}).
The benchmark solutions were obtained via front-tracking, using very fine discretizations for the parametrized curves.
Second-order convergence is evident, until the errors (measured using an approximate Hausdorff distance) become quite small. 

\begin{table}[H]
\begin{center}
\begin{tabular}{|c|c|l|c|l|c|}
\hline
$h$ & $nt$ & 4-petal Err. & Order & 6-petal Err. & Order\\
\hline
$\frac{1}{100}$ & $1$ & 0.0050554  & -- & 0.010446 & --\\
\hline
$\frac{1}{200}$ & $2$ & 0.0009674   & $2.4$ & 0.001862 & $2.5$\\
\hline
$\frac{1}{400}$ & $4$ & 0.00019335 & $2.3$ & 0.00033639 & $2.5$\\
\hline
$\frac{1}{800}$ & $8$ & $3.6891\times 10^{-5}$  & $2.4$ & $7.0263\times 10^{-5}$  & $2.3$\\
\hline
$\frac{1}{1600}$ & $16$ & $9.1215\times 10^{-6}$  & $2.0$ & $1.6355\times 10^{-5}$  & $2.1$\\
\hline
$\frac{1}{3200}$ & $32$ & $4.6133\times 10^{-6}$ & $1.0$ & $5.8529\times 10^{-6}$ & $1.5$\\
\hline
\end{tabular}
\end{center}
\caption{\footnotesize Convergence study for the four-petal and six-petal flower-shaped initial interface using the fast and accurate median computation described in Section \ref{sec:Fast} with kernel defined in (\ref{eq:2ndorderkernel}), exhibiting second order convergence in time until the error is down to $\approx 5 \times 10^{-6}$.}
\label{table:flowers}
\end{table}

An extension to three dimensions of the ideas explored in this and the previous section will be taken up in a subsequent study.

\section{An Obstruction}
\label{sec:barrier}
An immediate question is whether a second order in time median scheme for motion by mean curvature can be found in dimensions $d\geq 3$, in analogy with the two dimensional construction (\ref{eq:2ndorderkernel}) of Section \ref{sec:HighOrder}.
Here, we report some bad news on this front.
Section \ref{appendix:2} (Appendix 2) contains Taylor expansions (in the time step size $k$) for the exact solution of motion by mean curvature and one step of the threshold dynamics scheme (\ref{alg:td}) for an {\em arbitrary} radially symmetric kernel $K$ in dimension $d=3$.
They imply the following:

\begin{proposition}
There is no median filter scheme of the form (\ref{eq:obw}) that is second order consistent with motion by mean curvature equation (\ref{eq:levelset}) in dimensions $d\geq 3$.
\end{proposition}

\noindent {\bf Proof:}
We compare the Taylor expansion (\ref{eq:3dTaylor}) at short time $t$ for the exact solution of mean curvature motion of a surface given as the graph of a function, and that of the approximation generated after one time step of size $t$ with the median filter (\ref{eq:obw}).
The expansions are in powers of $t$, and describe where an initial interface passing through $(0,0,0)$ with unit normal $(0,0,1)$ intersects the $z$-axis at time $t$; see Appendix 2.

Notice that the various terms appearing in either expansion depend on the initial interface through the same geometric quantities: Only the quantities
\begin{equation}
H(0,0) \; , \; \Delta_S H(0,0) \; , \; H^3(0,0) \; \mbox{, and} \; H(0,0) \kappa(0,0)
\end{equation}
appear in either.
Matching the coefficients $B_0,\ldots,B_3$ defined in (\ref{eq:coeff}) that appear in the expansion (\ref{eq:boldh})
in Section \ref{appendix:2} for one step with the weighted median scheme to the corresponding coefficients in the Taylor expansion (\ref{eq:3dTaylor}) for the exact solution gives:
\begin{equation}
\label{eq:Bs}
    \begin{split}
        &B_0^2 = 2B_1\\
        &2B_1 = -B_2\\
        &6B_1 = B_3
    \end{split}
\end{equation}
The first equation of (\ref{eq:Bs}), $B_0^2 = 2 B_1$, implies 
\begin{equation}
\label{eq:solveeqn1}
    \left(\frac{\langle K \rangle_3 }{\langle K \rangle_1}\right)^2 = \frac{1}{2}\frac{\langle K \rangle_5 }{\langle K \rangle_1} 
\end{equation}
Plugging (\ref{eq:solveeqn1}) into the second equation of (\ref{eq:Bs}), $2B_1 = -B_2$, gives 
\begin{equation}
\begin{split}
       \frac{1}{32}\frac{\langle K \rangle_5 }{\langle K \rangle_1} &= \frac{5}{128}\frac{\langle K \rangle_5 }{\langle K \rangle_1} -\frac{1}{32} \left(\frac{\langle K \rangle_3 }{\langle K \rangle_1}\right)^2\\
       & = \frac{5}{128}\frac{\langle K \rangle_5 }{\langle K \rangle_1}  - \frac{1}{64} \frac{\langle K \rangle_5 }{\langle K \rangle_1}\\
       &= \frac{3}{128}\frac{\langle K \rangle_5 }{\langle K \rangle_1},
\end{split}
\end{equation}
which implies $\langle K \rangle_5 = 0$ and therefore $B_0 = 0$.
Hence, the interface remains stationary at the $O(1)$ time scale: not even first order consistency is achieved if one attempts to remove the leading order error terms.
$\Box$

\section{Multiphase Median Schemes}
\label{sec:Multiphase}

The threshold dynamics algorithm (\ref{eq:td1}) \& (\ref{eq:td2}) and its variational interpretation (\ref{eq:tdvariational}), extend to multiple phases \cite{merriman_bence_osher,mbo92,ruuth1,esedoglu_otto}, which is the motion by mean curvature of partitions.
Here, the domain $D\subset\mathbb{R}^d$ divided into at most $N$ essentially disjoint sets $\Sigma_1,\ldots,\Sigma_N$:
\begin{equation}
\label{eq:partition}
D = \bigcup_{j=1}^N \Sigma_j \mbox{ and } \Sigma_i \cap \Sigma_j = \partial\Sigma_i \cap \partial\Sigma_j \mbox{ for } i\not=j.
\end{equation}
The energy is of the form
\begin{equation}
\label{eq:mE}
E(\Sigma_1,\ldots,\Sigma_m) = \sum_{i,j} \sigma_{i,j} \mbox{Area}( \partial\Sigma_i \cap \partial\Sigma_j)
\end{equation}
with $\sigma_{i,i} = 0$ for $i\in\{1,2,\ldots,N\}$.
The parameters $\sigma_{i,j} = \sigma_{j,i}$ with $\sigma_{i,i} = 0$ are known as surface tensions; they are positive for $i\not= j$, and satisfy a triangle inequality:
\begin{equation}
\label{eq:triangle}
\sigma_{i,j} + \sigma_{j,k} \geq \sigma_{i,k} \mbox{ for any  distinct } i,j,k.
\end{equation}
The approximate energies become \cite{esedoglu_otto}:
\begin{equation}
\label{eq:amE}
E_K(\Sigma_1,\ldots,\Sigma_N) = \sum_{i,j} \sigma_{i,j} \int_{\Sigma_i} K * \mathbf{1}_{\Sigma_j} \, dx
\end{equation}
and lead to the following scheme:

\begin{algorithm}[H]
\caption{Multiphase Threshold Dynamics \cite{esedoglu_otto}}
\label{alg:mtd}
\begin{algorithmic}[1]
\STATE Convolution step:
        \begin{equation}
          \label{eq:mtd1}
          u_j^n = K_{k} * \left( \sum_{i=1}^N \sigma_{i,j} \mathbf{1}_{\Sigma_i^n} \right)
        \end{equation}
\STATE Thresholding step:
        \begin{equation}
          \label{eq:mtd2}
          \Sigma_i^{n+1} = \Big\{ x \, : \, u_i^n(x) = \min_j u_j^n(x) \Big\}.
        \end{equation}
\end{algorithmic}
\end{algorithm}

Step (\ref{eq:mtd2}) of Algorithm \ref{alg:mtd} can be augmented with a tie-break rule to ensure $\Sigma^{n+1}_i$ form a partition.
Let us write $\Sigma = (\Sigma_1,\Sigma_2,\ldots,\Sigma_N)$ and
\begin{equation}
\label{eq:vectS}
\mathbf{S}_K \Sigma^n  = \Sigma^{n+1}
\end{equation}
to denote one step of this scheme (it in fact agrees with definition (\ref{eq:T}) of $\mathbf{S}_K$ when $N=2$).
A large class of surface tensions $\sigma$ are identified in \cite{esedoglu_otto} that guarantee Algorithm \ref{alg:mtd} dissipates (\ref{eq:amE}).

In this section, we explore level-set versions of algorithm (\ref{eq:mtd1}) \& (\ref{eq:mtd2}); they can be viewed as multiphase analogues of the median filter.
The ideal situation would be to formulate a variational problem of the form (\ref{eq:mE}) at every level, and sum them up, much like in Section \ref{sec:Connection}.
However, if $\phi_1,\phi_2,\ldots,\phi_N$ are regular level-set functions representing the $N$ phases, i.e.
\begin{equation}
\mathbf{T}_0 \phi_i = \big\{ x \, : \, \phi_i(x) \geq 0 \big\} = \Sigma^0_i,
\end{equation}
then unfortunately other levels of $\phi_i$ do not form a partition: $\mathbf{T}_\lambda \phi_i = \{ x \, : \, \phi_i(x) \geq \lambda\}$ leave a vacuum when $\lambda > 0$, and overlap when $\lambda <0$.
Hence it is not clear how to extend the partition problem (\ref{eq:partition}) \& (\ref{eq:mE}) to all levels $\lambda$ so that the resulting analogue of scheme (\ref{eq:vectS}), call it $\widetilde{\mathbf{S}}_{K,\lambda}$, would also ensure a maximum principle of the form $\big( \widetilde{\mathbf{S}}_{K,\mu} \Sigma^n \big)_i \subseteq \big( \widetilde{\mathbf{S}}_{K,\lambda} \Sigma^n \big)_i$ for every $\mu \geq \lambda$ (so that the sets $\big( \widetilde{\mathbf{S}}_{K,\lambda} \Sigma^n \big)_i$ from different $\lambda$ but fixed $i$ can be ``stacked'' to form a function from which the sets can later be recovered via thresholding).
Perhaps vacuum can be treated as another phase, and overlaps penalized, in such a way that precludes their formation when initially absent.
We leave this line of inquiry to another time.

Here, we instead settle for a less satisfactory extension in which super level-sets other than $\lambda = 0$ do not appear to solve a specific variational problem that we are able to identify, but the following valuable properties hold:
\begin{itemize}
\item $\lambda=0$ super level-sets $\mathbf{T}_0 \phi_i$ still form a partition and evolve exactly by Algorithm \ref{alg:mtd},
\item A comparison principle for (reversibly) stacking super level-sets from different levels $\lambda$ holds for each individual phase.
\end{itemize}
We rely on the following simple partial comparison principle that holds for Algorithm \ref{alg:mtd}:
\begin{proposition}
\label{comparison}
Assume $\widetilde{\Sigma}_i^n \subseteq \Sigma_i^n$ and for all $j\not= i$, $\Sigma_j^n \subseteq \widetilde{\Sigma}_j^n$.
Then, under Algorithm \ref{alg:mtd}, $\widetilde{\Sigma}_i^{n+1} \subseteq \Sigma_i^{n+1}$.
\end{proposition}
\medskip
\noindent {\bf Proof:}
Let $u_\ell^n$ and $\widetilde{u}_\ell^n$ be defined as in step (\ref{eq:mtd1}) of Algorithm (\ref{alg:mtd}).
Also, let $\chi_\ell = K * \mathbf{1}_{\Sigma_\ell^n}$ and $\widetilde{\chi}_\ell = K*\mathbf{1}_{\widetilde{\Sigma}_\ell^n}$.
Noting $\sum_\ell \chi_\ell = \sum_\ell \widetilde{\chi}_\ell = 1$,  for $j\not = i$ we have:
\begin{equation}
\begin{split}
\widetilde{u}_i^n - \widetilde{u}_j^n =& u_i^n - u_j^n + \sum_{\ell \not= i} \sigma_{i,\ell} \Big( \widetilde{\chi}_\ell - \chi_\ell \Big) + \sum_{\ell \not = j} \sigma_{j,\ell} \Big( \widetilde{\chi}_\ell - \chi_\ell \Big)\\
=& u_i^n - u_j^n + \sum_{\ell \not= i,j} \big( \sigma_{i,\ell} + \sigma_{j,\ell} \big) \Big( \widetilde{\chi}_\ell - \chi_\ell \Big)\\
&+ \sigma_{i,j} \Big( \widetilde{\chi}_i - \chi_i \Big) + \sigma_{i,j} \sum_{\ell \not= j} \Big( \chi_\ell - \widetilde{\chi}_\ell \Big)\\
\geq& u_i^n - u_j^n + \sigma_{i,j} \sum_{\ell \not= i,j} \Big( \widetilde{\chi}_\ell - \chi_\ell \Big)\\
&+ \sigma_{i,j} \Big( \widetilde{\chi}_i - \chi_i \Big) + \sigma_{i,j} \sum_{\ell \not= j} \Big( \chi_\ell - \widetilde{\chi}_\ell \Big) \;\; \Big(\mbox{by (\ref{eq:triangle})} \Big)\\
=& u^n_i - u^n_j
\end{split}
\end{equation}
Step (\ref{eq:mtd2}) of Algorithm \ref{alg:mtd} now implies $u^{n+1}_i \leq u^n_i$. $\Box$
\medskip

Motivated by Proposition \ref{comparison}, our algorithm is as follows: Given $N$ level set functions $\phi(x) = \big( \phi_1(x),\ldots,\phi_N(x) \big)$ representing the initial phases $\Sigma^0_1,\ldots,\Sigma^0_N$ as $\Sigma^0_i = \{ x \, : \, \phi_i(x) \geq 0 \}$, we consider $N$ different partitions at every level $\lambda$: The partition at level $\lambda$ {\em from the point of view of phase $i$} is:
\begin{equation}
\mathcal{S}_i(\lambda) = \big\{ \Sigma_{i,1}(\lambda) \, , \, \Sigma_{i,2}(\lambda) \, , \, \ldots \, , \, \Sigma_{i,N}(\lambda) \big\}
\end{equation}
where
\begin{equation}
\Sigma_{i,j}(\lambda) =
\begin{cases}
\big\{ x \, : \, \phi_i(x) \geq \lambda \big\} & \mbox{ for } j=i,\\
\big\{ x \, : \, \phi_i(x) < \lambda \mbox{ and } \phi_j(x) \geq \max_{k\not\in\{i,j\}} \phi_k(x) \big\} & \mbox{ for } j\not= i.
\end{cases}
\end{equation}
Every super level-set $\mathbf{T}_\lambda\phi_i$ of phase $i$ will in effect be updated by applying one step of multiphase threshold dynamics (\ref{eq:mtd1}) \& (\ref{eq:mtd2}) to the partition $\mathcal{S}_i(\lambda)$.
To that end, let
\begin{equation}
\label{eq:psij}
\psi_{i,j} \phi(x,\lambda) = \sum_{k\not= j} \sigma_{j,k} \int_{\Sigma_{i,k}(\lambda)} K(x-y) \, dy.
\end{equation}
Then, by Proposition \ref{comparison}, the sets
\begin{equation}
\big\{ x \, : \, \psi_{i,i}\phi(x,\lambda) = \min_{j} \psi_{i,j} \phi(x,\lambda) \big\}
\end{equation}
are decreasing in $\lambda$.
We therefore let
\begin{equation}
\left( \mathbf{M}_K \phi \right)_i = \sup \Big\{ \lambda \, : \, \psi_{i,i} \phi(x,\lambda) = \min_{j} \psi_{i,j}\phi(x,\lambda) \Big\}.
\end{equation}
With these definitions, the multiphase algorithm is simply
\begin{equation}
\label{eq:m}
\phi^{n+1} = \mathbf{M}_K \phi^n.
\end{equation}
Then, we have 
\begin{equation}
\label{eq:tmst}
\mathbf{T}_0 \mathbf{M}_K \phi = \mathbf{S}_K \mathbf{T}_0 \phi
\end{equation}
i.e. the $0$-level sets are evolved precisely by multiphase threshold dynamics (\ref{eq:mtd1}) \& (\ref{eq:mtd2}).
Scheme (\ref{eq:m}) can be implemented using a sort operation, as follows:
\medskip

%\begin{algorithm}
%\caption{Multiphase Median Filter to update $\phi^n_i(x)$, $\sigma_{i,j} = \delta_{i,j}$:}
%\begin{algorithmic}[1]
%\STATE {\bf sort} the level set values $\{ \phi^n_i(y_1), \phi^n_i(y_2), \ldots, \phi^n_i(y_m)\}$ so that the permutation $p \, : \, \{1,2,\ldots,m\} \to \{1,2,\ldots,m\}$ satisfies
%$$ \phi^n_i \big( y_{p(1)} \big) \leq \phi^n_i \big( y_{p(2)} \big) \leq \cdots \leq \phi^n_i \big( y_{p(m)} \big).$$
%\STATE {\bf set} $C_j=0$ for $j=1,2,\ldots,N$ and $\ell = 1$.
%\WHILE{$\|K\|-C_i \geq \max_{j\not= i} C_j$}
%\STATE $C_i \leftarrow C_i + K \big( x + y_{p(\ell)} \big)$
%\STATE {\bf set} $\displaystyle j_* = \argmax_{j\not= i} \phi^n_j(y_{p(\ell)})$
%\STATE $C_{j_*} \leftarrow C_{j_*} + K \big( x + y_{p(\ell)} \big) $
%\STATE $\ell \leftarrow \ell+1$
%\ENDWHILE
%\STATE $ \phi^{n+1}_i(x) = \frac{1}{2} \Big( \phi^n_i(y_{p(\ell-1)}) + \phi^n_i(y_{p(\ell)}) \Big)$
%\end{algorithmic}
%\end{algorithm}

\begin{figure}
\begin{center}
\includegraphics[scale=0.3]{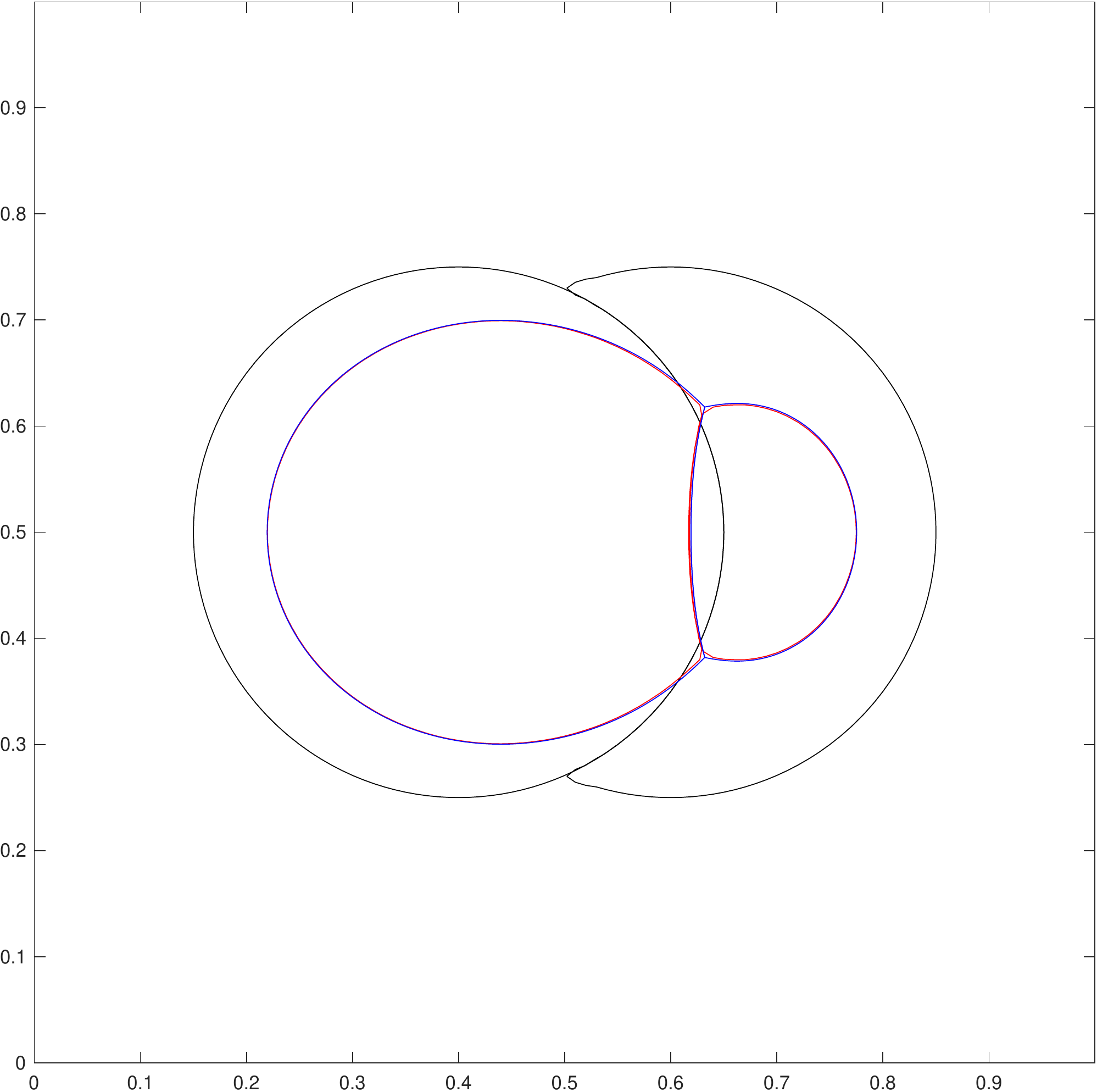}
\caption{\footnotesize Initial condition in black, multiphase median filter computed on $100\times 100$ grid with 40 time steps in red, and naive threshold dynamics on $4096\times 4096$ grid with 40 time steps in blue.}
\label{fig:3phase1}
\end{center}
\end{figure}

\begin{figure}
\begin{center}
\includegraphics[scale=0.15]{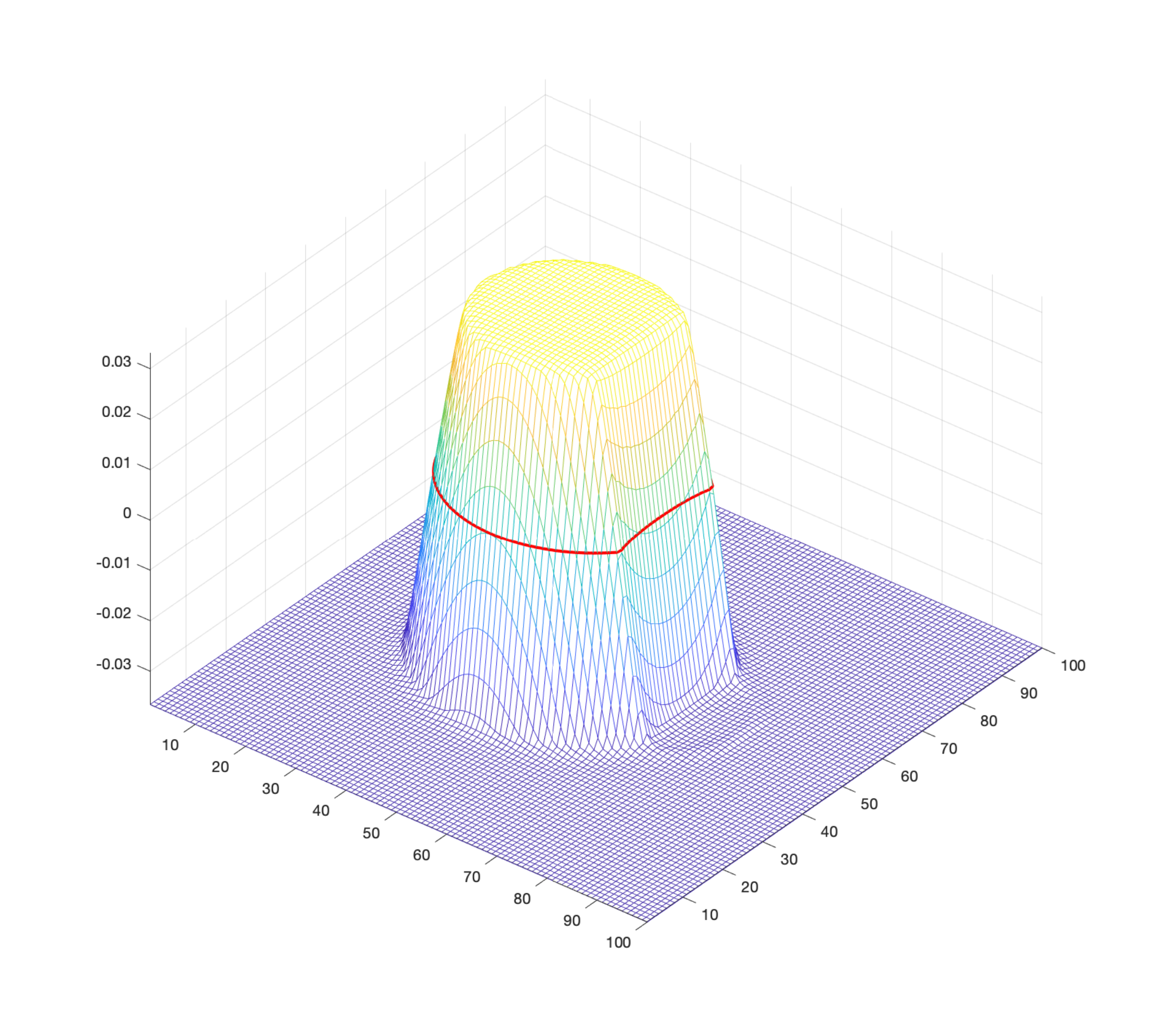}
\includegraphics[scale=0.15]{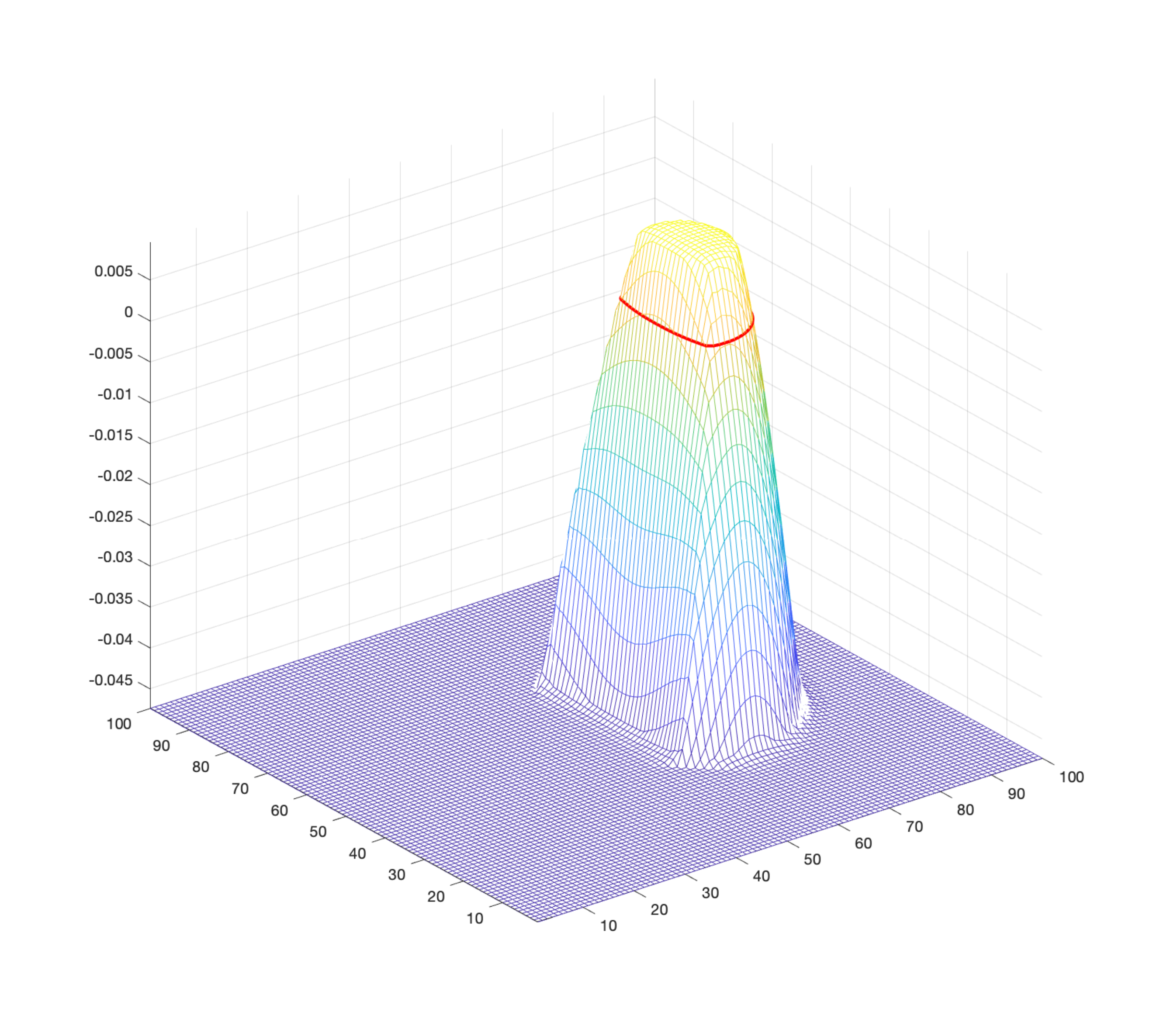}
\includegraphics[scale=0.15]{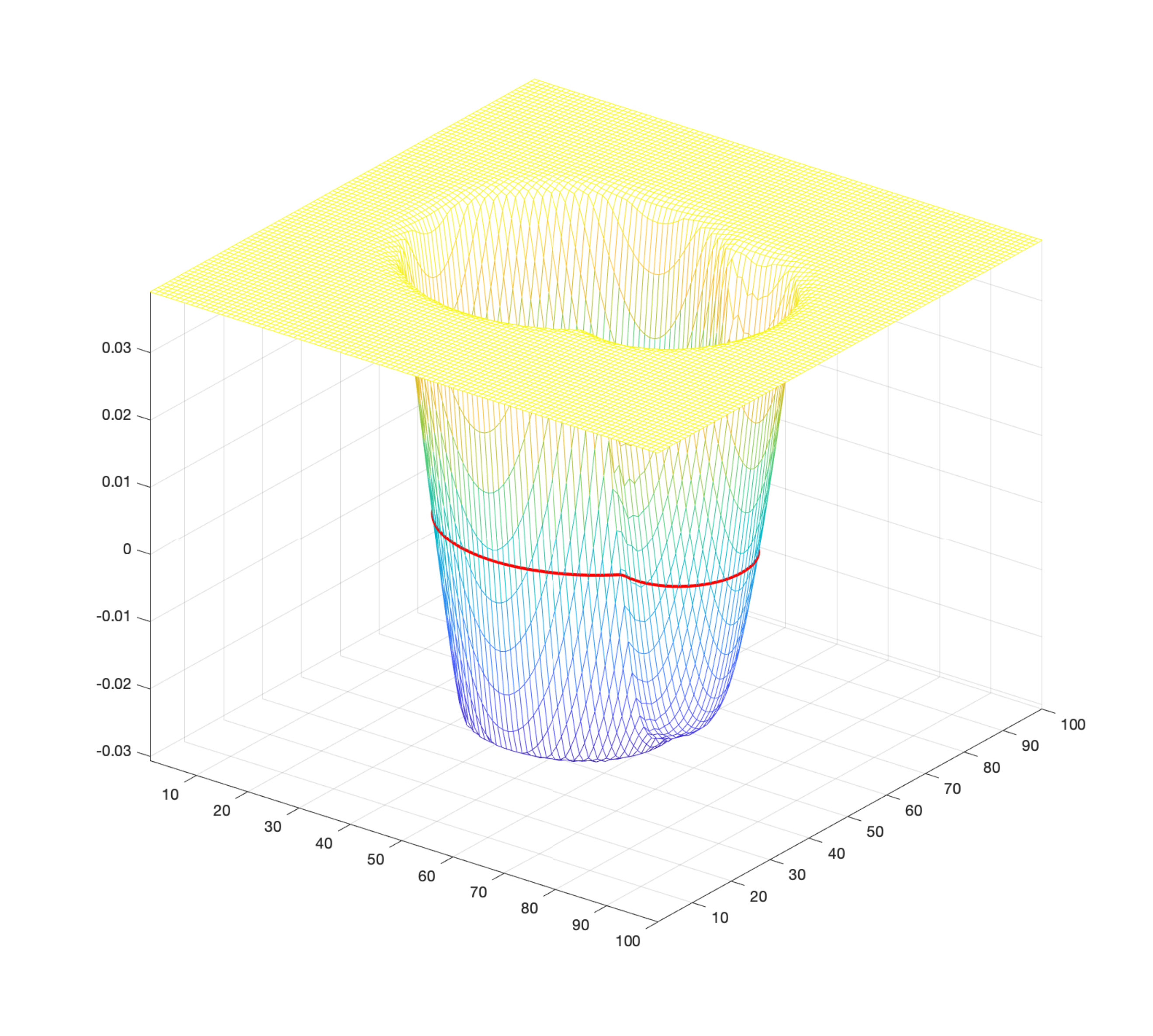}
\caption{\footnotesize Level set functions for the multiphase median filter computation in Figure \ref{fig:3phase1} at final time, with narrow banding but no reinitialization.
They appear to remain regular ($|\nabla \phi_j|$ not too small or large) at least for a while near the relevant $0$-level set, shown in red.}
\label{fig:3phase2}
\end{center}
\end{figure}

\begin{figure}[H]
\begin{center}
\includegraphics[scale=0.415]{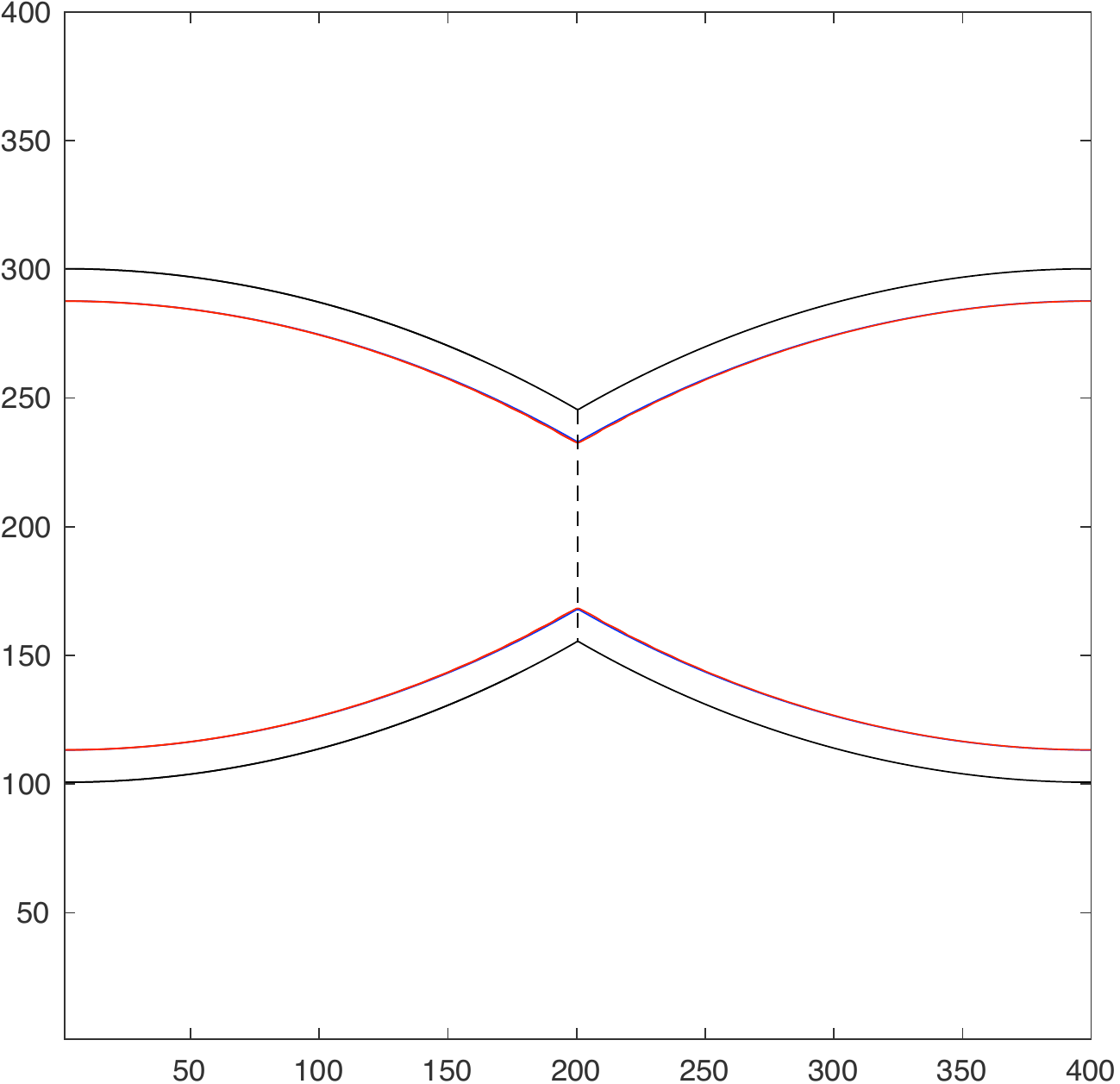}
\hspace{0.5 cm}
\includegraphics[scale=0.35]{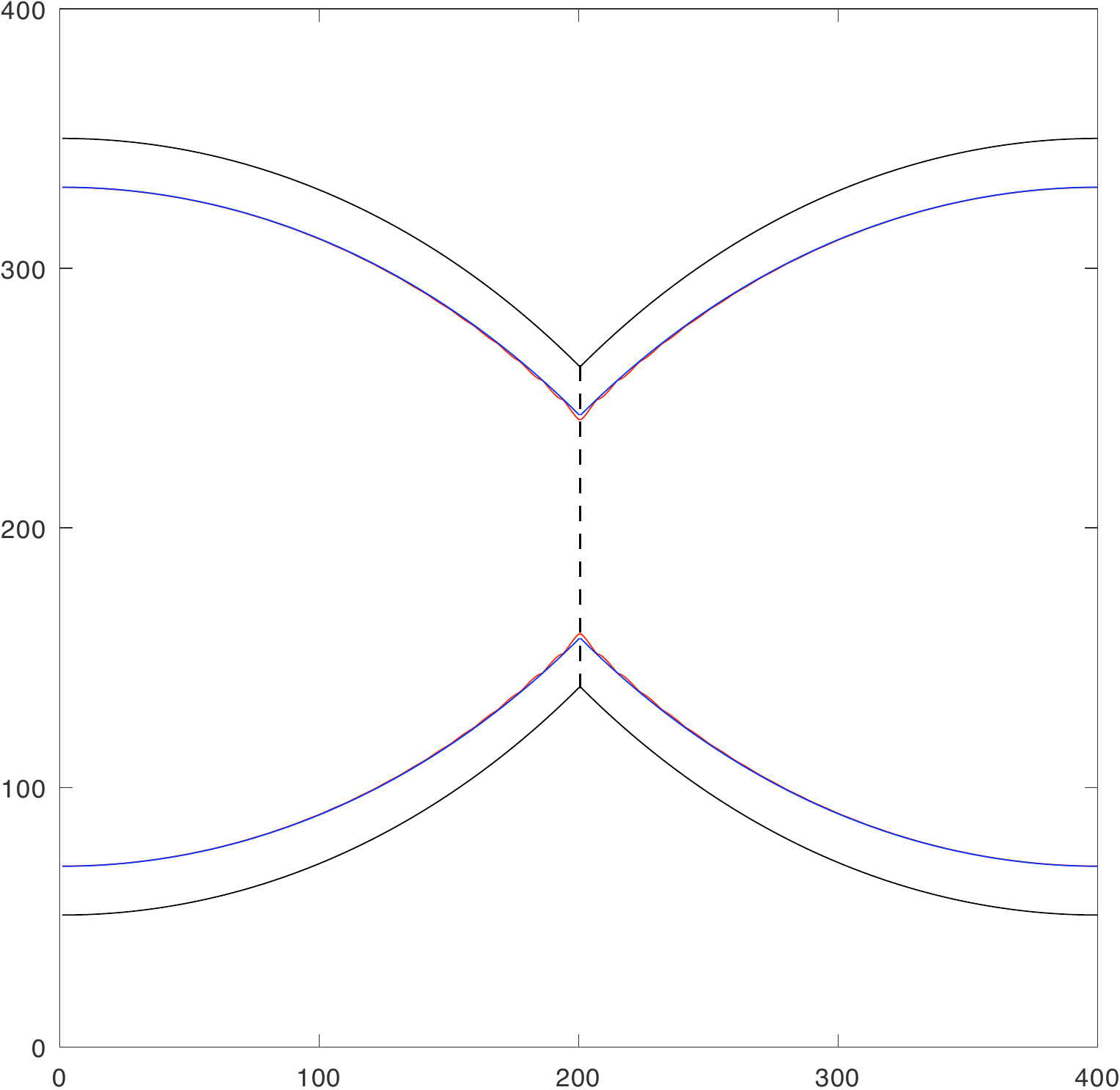}
\caption{\footnotesize Convergence study for the multiphase median filter on the family of exact solutions for three phase curvature motion known as ``grim reapers''.
Black curves are the initial condition. Blue are the exact solution at time $T=0.03$, and red are the approximate solution computed by multiphase median filter.
Left: Surface tensions chosen to give $(120^\circ,120^\circ,120^\circ)$ angle condition at the triple junction. Right: Surface tensions imply $(135^\circ,135^\circ,90^\circ)$ angle condition.}
\label{fig:grim}
\end{center}
\end{figure}

\begin{algorithm}[H]
\caption{A Multiphase Median Filter}
\label{alg:mmf}
Update $\phi^n_i(x)$ as follows:
\begin{algorithmic}[1]
\STATE {\bf sort} the level set values $\{ \phi^n_i(x+y_1), \phi^n_i(x+y_2), \ldots, \phi^n_i(x+y_m)\}$ so that the permutation $p \, : \, \{1,2,\ldots,m\} \to \{1,2,\ldots,m\}$ satisfies
$$ \phi^n_i \big( x+y_{p(1)} \big) \leq \phi^n_i \big( x+y_{p(2)} \big) \leq \cdots \leq \phi^n_i \big( x+y_{p(m)} \big).$$
\STATE {\bf set} $C_j = \sigma_{i,j} \| K \| $ for $j=1,2,\ldots,N$ and $\ell = 1$.
\WHILE{$C_i \leq \min_{j\not= i} C_j$}
\STATE {\bf set} $\displaystyle j_* = \argmax_{j\not= i} \phi^n_j(x+y_{p(\ell)})$
\STATE $C_j \leftarrow C_j + \big( \sigma_{j,j_*} - \sigma_{i,j} \big) K \big( y_{p(\ell)} \big) $ for each $j=1,2,\ldots,N$.
\STATE $\ell \leftarrow \ell+1$
\ENDWHILE
\STATE $ \phi^{n+1}_i(x) = \frac{1}{2} \Big( \phi^n_i(x+y_{p(\ell-1)}) + \phi^n_i(x+y_{p(\ell)}) \Big)$
\end{algorithmic}
\end{algorithm}
Just like Algorithms \ref{alg:ob} and \ref{alg:wmf}, Algorithm \ref{alg:mmf} thus applies threshold dynamics (in this case, the multiphase version given in Algorithm \ref{alg:mtd}) to infinitely many partitions at the same time.
In practice, good accuracy requires a large number $m$ of points $y_1,\ldots,y_m$ sampling the support of the kernel $K$; in numerical experiments, we therefore adopt the approach of Section \ref{sec:Fast} to approximate the integrals (\ref{eq:psij}) with weight $K=\mathbf{1}_{B_r(0)}$. 
Figure \ref{fig:3phase1} shows a three-phase computation where surface tensions are $\sigma_{i,j} = \delta_{i,j}$, corresponding to symmetric, $(120^\circ,120^\circ,120^\circ)$ junctions.
Figure \ref{fig:3phase2} shows the level set functions at final time; they appear to remain fairly regular, with their gradients $|\nabla\phi_i|$ neither too large nor too small near the $0$-level set, without resorting to reinitialization \cite{osher_sethian}, which can move the interface \cite{russo_smereka}.
Figure \ref{fig:grim} and Table \ref{table:grim} show a convergence study on a class of three-phase exact solutions known as ``grim reapers''.
These feature a curved interface of the form $y=f(x,t)$ on $x\in[0,\frac{1}{2}]$ with
\begin{equation}
\label{eq:exactgrim}
f(x,t) = \frac{\alpha}{\pi} \log \left( \cos \left( \frac{\pi}{\alpha}x \right) \right) - \frac{\pi}{\alpha} t
\end{equation}
that translates in the $y$-direction with constant speed $\frac{\pi}{\alpha}$.
As such, the junction at $x=\frac{1}{2}$ has angles $$\left( \frac{\pi(\alpha-1)}{\alpha},\frac{\pi(\alpha+1)}{2\alpha},\frac{\pi(\alpha+1)}{2\alpha} \right).$$
The experiments are with $\alpha=3$ corresponding to a $(120^\circ,120^\circ,120^\circ)$ junction, and with $\alpha=2$ corresponding to a $(90^\circ,135^\circ,135^\circ)$ junction the surface tensions of which are $\sigma_{1,3}= \sigma_{2,3} = \frac{\sqrt{2}}{2}$ and $\sigma_{1,2} = 1$.
The expected order of convergence as time step size vanishes, same as that of threshold dynamics in the presence of junctions, is $\frac{1}{2}$.
The experiments, which refine the spatial grid and time step size concurrently, exhibit slightly better rates; otherwise, they provide good evidence of convergence and reach errors that are quite small, measured as the $L^2$ norm of the difference between the computed and exact profile (\ref{eq:exactgrim}) at final time $T = 0.03$.

\begin{table}[H]
\begin{center}
\begin{tabular}{|c|c|c|c|c|c|}
\hline
$h$ & $nt$ & $(120,120,120)^\circ$ Err & Order & $(90,135,135)^\circ$ Err & Order\\
\hline
$1/50$ & $5$ & $0.00295$ & $--$ & $0.00441$ & $--$\\
\hline
$1/100$ & $10$ & $0.00187$ & $0.66$ & $0.00265$ & $0.74$\\
\hline
$1/200$ & $20$ & $0.00119$ & $0.65$ & $0.00163$ & $0.70$\\
\hline
$1/400$ & $40$ & $0.000727$ & $0.72$ & $0.00106$ & $0.62$\\
\hline
$1/800$ & $80$ & $0.000481$ & $0.59$ & $0.00053$ & $1$\\
\hline
\end{tabular}
\caption{\footnotesize Convergence study on ``grim reapers'', a family of exact solutions to three-phase motion by mean curvature, shown in Figure \ref{fig:grim}.
The order of convergence expected asymptotically is $\frac{1}{2}$.}
\label{table:grim}
\end{center}
\end{table}

\section{Appendix 1}
\label{sec:Appendix}
We repeat the calculation in \cite{esedoglu_guo} with the convolution kernel given by a delta function concentrated along a circle (instead of a Gaussian) in two dimensions, and then also consider a spherical shell in three dimensions.
Let the interface be given as the graph of the function $f(x)$, where $f(0)=0$ and $f'(0)=0$.
Recall that in $d=2$, curvature motion of an interface described as the graph of a function $u(x,t)$ is given by the PDE
\begin{equation}
\begin{split}
u_t &= \frac{u_{xx}}{1+u_x^2}\\
u(x,0) &= f(x).
\end{split}
\end{equation}
Taylor expanding $u(0,t)$ in $t$ at $t=0$ gives
\begin{equation}
\label{eq:exacttaylor}
u(0,t) = t f''(0) + t^2 \left( \frac{1}{2} f^{(iv)}(0) - \big( f''(0)\big)^3(0) \right) + O(t^3).
\end{equation}

Now consider the intersection of the graph with the circle $X^2 + (Y-h)^2 = r^2$, with $h=O(r^2)$ as $r\to 0$.
\begin{figure}[h]
\begin{center}
\includegraphics[scale=0.5]{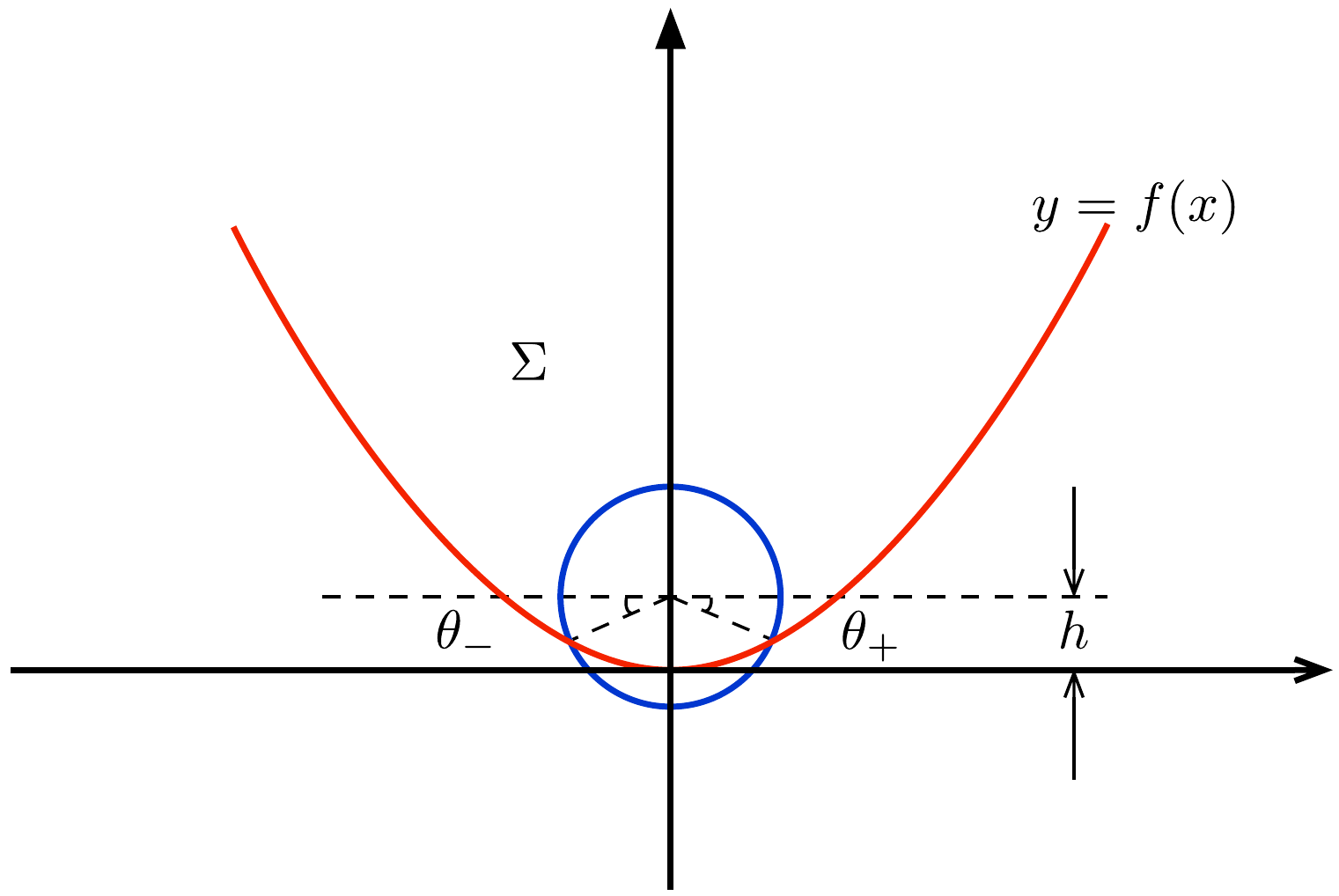}
\end{center}
\end{figure}
To find their intersection, solve the two equations
\begin{equation}
\begin{split}
X_\pm &= \pm r \sqrt{ 1 - \frac{(Y-h)^2}{r^2}} = \pm \left\{ 1 - \frac{(Y-h)^2}{2r^2} - \frac{(Y-h)^4}{8r^2} + \cdots \right\} \\
Y_\pm &= \frac{1}{2} f''(0) X^2 + \frac{1}{6} f'''(0) X^3 + \cdots
\end{split}
\end{equation}
by first eliminating $X$ in the second equation using the first, and then solving for $Y$ in the simplified second equation using the ansatz
\begin{equation}
Y_\pm = c_{\pm 2} r^2 + c_{\pm 3} r^3 + c_{\pm 4} r^4 + O(r^5)
\end{equation}
After the substitution $h=\mathbf{h} r^2$, one gets
\begin{multline}
\label{eq:Yintersection}
Y_\pm = \frac{1}{2} f''(0) r^2 \pm \frac{1}{6} f'''(0) r^3 + \left( -\frac{1}{8} \big( f''(0) \big)^3 \right.\\
\left. + \frac{1}{2} \big( f''(0) \big)^2 \mathbf{h} - \frac{1}{2} f''(0) \mathbf{h}^2 + \frac{1}{24} f^{(iv)}(0) \right) r^4 + O(r^5)
\end{multline}
Letting $\theta_\pm = \arcsin((Y_\pm-h)/r)$, this gives the expansion
\begin{equation}
\label{eq:thetatot}
\begin{split}
\theta_+ + \theta_- =& \left(-2 \mathbf{h}+f''(0)\right) r + \left(-\frac{1}{3} \mathbf{h}^{3}-\frac{1}{2} f''(0) \,\mathbf{h}^{2}+\frac{3}{4} \big(f''(0)\big)^2 \mathbf{h} \right.\\
& \left. -\frac{5}{24} \big(f''(0)\big)^{3}+\frac{1}{12} f^{(iv)}(0) \right) r^{3}
\end{split}
\end{equation}

\section{Appendix 2}
\label{appendix:2}
In this section, we extend the consistency calculation (Taylor expansions) from the previous Appendix to three dimensions.
In particular, we first obtain the expansion in case the convolution kernel is singular, concentrated on a sphere (in analogy with the circular kernel in two dimensions), and then treat the general case of a radially symmetric kernel.

Mean curvature motion of an interface described as the graph of a function $u(x,y,t)$ is given by the PDE:
\begin{equation}
\begin{split}
&u_t = \frac{ (1+u_y^2)u_{xx} + (1+u_x^2) u_{yy} - 2u_x u_y u_{x,y}}{1 + u_x^2 + u_y^2}\\
&u(x,y,0) = g(x,y).
\end{split}
\end{equation}
Taylor expanding $u(0,0,t)$ in $t$ at $t=0$ gives

\begin{equation}
 \label{eq:3dTaylor}
\begin{split}
u(0,0,t)  &=t H(0,0) \\
&\;\;\;+ \frac{1}{2} t^2 \left\{ \Delta^2 g(0,0) - 2H^3(0,0) + 6H(0,0) \kappa(0,0) \right\} + O(t^3)\\
&= t H(0,0) \\
&\;\;\;+ \frac{1}{2} t^2 \left\{  \Delta_S H(0,0) + H^3(0,0) - 2H(0,0) \kappa(0,0) \right\} + O(t^3)
\end{split}
\end{equation}

\noindent where
\begin{equation}
\label{eq:HKSLH}
\begin{split}
H(0,0) &= \Delta g(0,0) = g_{xx}(0,0) + g_{yy}(0,0)\\
\kappa(0,0) &= g_{xx}(0,0) g_{yy}(0,0) - g^2_{xy}(0,0)\\
\Delta_S H(0,0) &= \Delta^2 g(0,0) - 3H^3(0,0) + 8 H(0,0) \kappa(0,0).
\end{split}
\end{equation}
$H(0,0)$ and $\kappa(0,0)$ are the mean curvature and Gaussian curvature at the origin, respectively, and $\Delta_S H(0,0)$ is the surface Laplacian. 
We first consider using a spherical shell with a radius $r$ as the convolution kernel. Denote this spherical shell by $S_r$. Then the  analogous equation  to (\ref{eq:conv2d}) used to solve for $\mu$ becomes:
\begin{equation}
\label{eq:sphericalconv}
    \Big(\mathbf{1}_{\{(x,y,z) \, : \, z \geq g(x,y) \}} * S_r \Big)(0,0,\mathbf{h} r^2) = 2\pi r^2
\end{equation}
The left-hand side convolution in (\ref{eq:sphericalconv}) is the area of the spherical shell surface that lies in the interior of the hyperplane:
\begin{equation}
\label{eq:spherical_shell}
\begin{split}
      A(r) &:=  r^2 \int_0^{2\pi} \int_{0}^{\frac{\pi}{2}-\theta_+(\Theta)} \sin(\theta) d\theta d\Theta\\
      &= r^2\int_0^{2\pi} 1 - \cos\left(\frac{\pi}{2}-\theta_+(\Theta)\right) d\Theta\\
      &= r^2\int_0^{2\pi} 1 - \sin\left(\theta_+(\Theta)\right) d\Theta
\end{split}
\end{equation}
%\begin{figure}[h]
%\begin{center}
%\includegraphics[scale=0.5]{figs/3d_int2.pdf}
%\end{center}
%\end{figure}
where $\theta$ and $\Theta$  denote the azimuthal angle and the polar angle in the spherical coordinate system, respectively.
We can rotate the interface $z = g(x,y)$ clockwise by $\Theta$ ($0\leq\Theta<2\pi$), and only repeat the calculation done in Appendix \ref{sec:Appendix} in the y-z plane.
Denote the intersection between the rotated interface and the y-z plane by $z = G^{\Theta}(y)$, then the derivatives of $G^{\Theta}$ are given by:
%where $\cos(\phi_+(\theta))$ can be found exactly using (\ref{eq:thetatot}), which is $\frac{z_+ - h}{r}$. Then rotate $g(x,y)$ by $\theta$ to obtain the expression for $1 - \cos\left(\phi_+(\theta)\right)$ for each $\theta$. Without loss of generality, we can repeat the calculation done in Appendix \ref{sec:Appendix} by rotating every point in rotate each point lies in the graph of the function $g(x,y)$ to the $y-z$ plane, located in the function $z = G(y)$, then 
\begin{equation}
    \begin{split}
        G_{yy}^{\Theta} &= \sin(\Theta)^2 g_{xx} - 2\sin(\Theta)\cos(\Theta) g_{xy} +\cos(\Theta)^2 g_{yy}\\
        G_{yyy}^{\Theta} &= -\sin(\Theta)^3 g_{xxx} + 3\sin(\Theta)^2\cos(\Theta) g_{xxy} \\
        &\qquad- 3\sin(\Theta)\cos(\Theta)^2 g_{xyy} + \cos(\Theta)^3 g_{yyy}\\
        G_{yyyy}^{\Theta} &= \sin(\Theta)^4 g_{xxxx} - 4\sin(\Theta)^3\cos(\Theta) g_{xxxy} \\
        &\qquad+ 6\sin(\Theta)^2\cos(\Theta)^2 g_{xxyy} -4\sin(\Theta)\cos(\Theta) ^3 g_{xyyy} \\
        &\qquad+ \cos(\Theta)^4 g_{yyyy}
    \end{split}
\end{equation}

Plug the expressions of the derivatives into the adapted version of (\ref{eq:Yintersection}):
\begin{equation}
    \begin{split}
        1-\sin\left(\theta_+(\Theta)\right)  = &1 - \left(\frac{G_{yy}^{\Theta}}{2} - \mathbf{h}\right)r - \frac{G_{yyy}^{\Theta}}{6}r^2 -\\
        & \left(-\frac{1}{2}\mathbf{h}^2G_{yy}^{\Theta}+\frac{1}{2}\mathbf{h}(G_{yy}^{\Theta})^2 -\frac{1}{8}(G_{yy}^{\Theta})^3 +\frac{1}{24}G_{yyyy}^{\Theta}\right)r^3
    \end{split}
\end{equation}

Write 
\begin{equation}
\label{eq:defboldh}
\mathbf{h} = k_0 + k_1 r + k_2 r^2 + k_3 r^3 + \ldots
\end{equation}
and integrate over $\Theta$; we find:
\begin{equation}
\begin{split}
        A(r) &= 2\pi r^2 + r^3\left( \frac{1}{2}\pi H(0,0) - 2 \pi k_0\right)\\
        & + 2r^4  k_1\pi\\
        & + r^5 \pi\left( 2 k_2 - \frac{1}{32}\Delta^2 g(0,0) + \frac{5}{64} H^3(0,0) -\frac{1}{16} H(0,0)\kappa(0,0)\right.\\
        &\qquad\left.-\frac{3}{8}k_0 H(0,0)^2 + \frac{k_0}{2} \kappa(0,0) + \frac{k_0^2}{2}H(0,0)\right)\\
        & + h.o.t.\vphantom{r^2}
\end{split}
\end{equation}

Set $A(r) - 2\pi r^2 = 0$, we find
\begin{equation}
    \begin{split}
        k_0 &= \frac{1}{4}H(0,0),\\
        k_1 &= 0\\
        k_2 &= \frac{1}{64} \Delta^2 g(0,0) - \frac{1}{128} H^3(0,0) + \frac{1}{32} H(0,0)\kappa(0,0)\\
        &= \frac{1}{64}\Delta_S H(0,0) +\frac{5}{128} H^3(0,0) - \frac{3}{32}H(0,0)\kappa(0,0)
    \end{split}
\end{equation}
Plugging these into (\ref{eq:defboldh}) gives
\begin{equation}
\label{eq:boldhsphere}
\mathbf{h} = \frac{1}{4}H(0,0) + r^2 \Big( \frac{1}{64} \Delta^2 g(0,0) - \frac{1}{128} H^3(0,0) + \frac{1}{32} H(0,0)\kappa(0,0) \Big) + h.o.t.
\end{equation}
for the location of the interface along the $z$-axis after one time step with Algorithm \ref{alg:td} using a convolution kernel $K$ concentrated on a spherical shell; this special case is of interest for the original median filter scheme (Algorithm \ref{alg:ob}) from \cite{oberman2004} in dimension $d=3$.

We now generalize this expansion to any arbitrary radial symmetric kernel. The convolution becomes:
\begin{equation}
\label{eq:3dconvolution}
    \int_{0}^{\infty} K\left(\frac{r'}{\eps}\right) A(r')dr'.
\end{equation}
Here we divide radius $r'$ by $\eps (= O(\sqrt{dt}))$ since the kernel is expected to concentrate around $0$ as $dt\rightarrow 0$. 
Then introduce a new variable $r = \frac{r'}{\eps}$, the integral becomes
\begin{equation}
\label{eq:3dint}
     \int_{0}^{\infty} \eps K(r) A(\eps r) dr
\end{equation}

\noindent Define 
\begin{equation}
    \langle K \rangle_n : = \int_{0}^{\infty} r^n K(r) dr
\end{equation}
The integral can be expressed as:
\begin{equation}
    \begin{split}
        (\ref{eq:3dint}) = &2\pi \eps^3 \langle K \rangle_2 \\
        & + \eps^4\left(\frac{-2\pi H \langle K \rangle_3 + 8\pi \langle K \rangle_1 k_0}{4}\right) \\
        & + 2 \pi  \eps^5 k_1 \langle K \rangle_1\\
        & + \eps^6 \pi \left\{ \left( -\frac{1}{32} \Delta^2 g(0,0) + \frac{5}{64}H^3(0,0) - \frac{3}{16}H(0,0)\kappa(0,0)\right) \langle K \rangle_5\right. \\
        & +  \left(-\frac{3}{8}H^2(0,0)+\frac{1}{2}\kappa(0,0)\right) k_0\langle K \rangle_3\\
        &\left. + \left(\frac{1}{2}H(0,0) k_0^2 + 2 k_2\right)\langle K \rangle_1 \right\}
    \end{split}
\end{equation}
See \cite{grzibovskis_heintz,heintz05} who also obtain this expansion for the convolution (\ref{eq:3dconvolution}), but not the expansion (\ref{eq:boldh}) below for the location of the interface after one time step up to the order that we require.
(Also see e.g. \cite{mascarenhas,ruuth0} for earliest contributions, for $K$ given by a Gaussian).
Set the coefficients of $\eps^4$,  $\eps^5$,  $\eps^6$ to be zero, we have:

\begin{equation}
    \begin{split}
    k_0 &= \frac{\langle K \rangle_3 }{4\langle K \rangle_1 } H(0,0)\\    
    k_1 &= 0\\
    k_2& = \frac{1}{64}\frac{\langle K \rangle_5 }{\langle K \rangle_1 } \Delta^2 g(0,0)\\
    & -\left(\frac{5}{128}\frac{\langle K \rangle_5 }{\langle K \rangle_1 } - \frac{1}{32}\frac{\langle K \rangle_3^2 }{\langle K \rangle_1^2}\right) H^3(0,0)\\
    & + \left(\frac{3}{32}\frac{\langle K \rangle_5 }{\langle K \rangle_1} - \frac{1}{16} \frac{\langle K \rangle_3^2 }{\langle K \rangle_1^2}\right)H(0,0)\kappa(0,0)
    \end{split}
\end{equation}

We define:  
\begin{equation}
\label{eq:coeff}
    \begin{split}
        &B_0 := \frac{\langle K \rangle_3 }{4\langle K \rangle_1 },\\
        &B_1 := \frac{1}{64}\frac{\langle K \rangle_5 }{\langle K \rangle_1 },\\
        &B_2 := - \left(\frac{5}{128}\frac{\langle K \rangle_5 }{\langle K \rangle_1 } - \frac{1}{32}\frac{\langle K \rangle_3^2 }{\langle K \rangle_1^2}\right), \\
        &B_3 := \frac{3}{32}\frac{\langle K \rangle_5 }{\langle K \rangle_1} - \frac{1}{16} \frac{\langle K \rangle_3^2 }{\langle K \rangle_1^2}  
    \end{split}
\end{equation}
and therefore the expansion for $\mathbf{h}$ is:
\begin{equation}
\label{eq:boldh}
\begin{split}
     \mathbf{h} =& B_0 H(0,0) + r^2 \left(B_1 \Delta^2 g(0,0) + B_2 H(0,0)^3 + B_3 H(0,0)\kappa(0,0)\right) \\
     &+ h.o.t. 
\end{split}
\end{equation}

\section{Acknowledgements}
Selim Esedo\=glu and Jiajia Guo were supported by NSF DMS-2012015.
David Li was supported by the REU (research experience for undergraduates) component of the same grant.

\bibliographystyle{plain}
\bibliography{references}

\end{document}